\documentclass{amsart}

\usepackage{amsthm,amsfonts,amsmath,amssymb,latexsym,epsfig}
\usepackage{upref,amssymb,eucal}

\newtheorem{theorem}{Theorem}
\newtheorem{lemma}[theorem]{Lemma}
\newtheorem{proposition}[theorem]{Proposition}
\newtheorem{corollary}[theorem]{Corollary}
\newenvironment{example}{\medskip \refstepcounter{theorem}
\noindent  {\bf Example \thetheorem}.\rm}{\,}
\newenvironment{remark}{\medskip \refstepcounter{theorem}
\noindent  {\bf Remark \thetheorem}.\rm}{\,}
\newenvironment{definition}{\medskip \refstepcounter{theorem}
\noindent  {\bf Definition \thetheorem}.\rm}{\,}

\renewcommand{\thetheorem}{\thesection.\arabic{theorem}}

\theoremstyle{definition}

\theoremstyle{remark}

\def \cJ{\mathcal J}

\def \cN{{\mathcal N}}

\def \whk{^{\widehat{p_k}}}
\def \w1{^{\widehat{p_1}}}
\def \w2{^{\widehat{p_2}}}

\def \bZ{{\mathbb Z}}
\def \bG{{\mathbb G}}

\def \12{\delta_{\cP}}

\def \<{\langle}
\def \>{\rangle}

\def \cP{\mathcal P}

\def \cN{\mathcal N}

\def \cG{\mathcal G}

\def \bZ{{\mathbb Z}}

\def \bF{{\mathbb F}}

\def \bC{{\mathbb C}}

\def \bA{{\mathbb A}}

\def \cO{\mathcal O}
\def \ra{\rightarrow}

\def \bX{{\bf X}}

\def \bF{{\mathbb F}}

\def \bQ{{\mathbb Q}}

\def \bQ{{\mathbb Q}}

\def \bZ{{\mathbb Z}}
\def \bG{{\mathbb G}}

\begin{document}

\title{Arithmetic Laplacians}
\author{Alexandru Buium}
\thanks{During the preparation of this work, the first author
was partially supported by NSF grant DMS 0552314.}
\author{Santiago R. Simanca}
\address{Department of Mathematics and Statistics \\
University of New Mexico \\ Albuquerque, NM 87131}
\email{buium@math.unm.edu, santiago@math.unm.edu}

\begin{abstract}
We develop  an arithmetic analogue of elliptic partial differential equations.
The r\^ole of the space coordinates is played by a family of primes, and
that of the space derivatives along the various primes are  played by
corresponding Fermat quotient operators subjected to certain commutation
relations. This
leads to arithmetic linear partial differential equations on algebraic
groups that are analogues of certain operators in analysis constructed
from Laplacians. We classify all such equations on one dimensional groups,
and analyze their spaces of solutions.
\end{abstract}

\maketitle

\section{Introduction}
\subsection{Main concepts and results}
In a series of articles (beginning with \cite{char}),
the first author developed an arithmetic analogue of ordinary differential
equations (ODEs); cf. \cite{book} for an account of this theory.
The aim of the present work is to extend some of this theory to the
partial differential case by considering arithmetic analogues of
elliptic partial differential equations (PDEs).
A different extension of the ODE  theory, in the direction of parabolic and
hyperbolic PDEs, was developed by the authors in \cite{pde, pdmod}. The
work here relies on the ODE theory  but is independent of \cite{pde,pdmod}.

Before explaining the main results of this paper let us put our present work
in perspective. The general aim of this paper and of \cite{pde,pdmod}
is to pass from one  variable to several variables.
The independent variables we have in mind are of two
types, the ``geometric independent variables,'' whose number we denote by
an integer $d_1\geq 0$, and the ``arithmetic independent variables,''
whose number we denote by the integer $d_2\geq 1$. We would like to extend
the ODE theory to a PDE theory in $d_1+d_2$ independent variables.
We explain this in some detail next. For convenience, we introduce a third
integer parameter $d_3\geq 1$ to denote the number of ``dependent variables''
that may arise in pursuing the said extension.

The ODE theory in \cite{char, book} corresponds to the case of $0+1$
independent variables, where the one independent arithmetic variable is
represented by a prime integer $p$. The role of the derivative with respect
to $p$
is played by a Fermat quotient operator $\delta_p$ that on integers
$a \in \bZ$ acts as
$$
\delta_p a=\frac{a-a^p}{p}\, .
$$
If we view this theory in a manner analogous to particle mechanics,
as in \cite{book}, then the prime $p$ is the analogue of ``time,'' and
$d_3$ is the dimension of the ``configuration space.'' If we pursue the
analogy with field theory instead, then $p$ is the analogue of a
``space'' variable and $d_3$ is the dimension of the
``space of internal states.''

The PDE theory in \cite{pde, pdmod} corresponds to the case of $1+1$
independent variables. The one geometric variable, $q$, is viewed as the
exponential of ``time,''  and the one arithmetic variable, viewed as the
analogue of space, is represented, again, by a prime $p$. The
derivative in the geometric direction is given by the usual derivation
$\delta_q=qd/dq$, while the role of the derivative with respect to $p$ is
again played by a Fermat quotient operator $\delta_p$.
The main purpose of \cite{pde} was to determine
all linear PDEs on the algebraic groups of dimension $d_3=1$
(precisely, on the group schemes $\bG_a$, $\bG_m$, $E$, where $E$ is any
elliptic curve), and to study the solutions of all such arithmetic linear
PDEs. In line with \cite{char,book}, these arithmetic linear PDEs are
referred to as $\{\delta_p,\delta_q\}$-{\it characters}. They can be viewed
as arithmetic analogues of evolution equations in analysis; for $\bG_m$,
the basic examples are analogues of convection equations
whereas for elliptic curves, the basic examples are analogues of
convection, heat, and wave equations. In \cite{pdmod}, the one dimensional
groups were replaced by modular curves. The resulting equations now can be
viewed as arithmetic analogues of convection and heat equations.
(Wave equations do not have analogues in the modular case.)

In this article, we develop an arithmetic PDE theory in $0+d_2$
independent variables, where $d_2 \geq 2$, integer which for simplicity
we now simply denote by $d$.
We view the new equations as independent of the geometric time, but the
analogue of the space coordinates is now a finite set
$\cP=\{p_1,\ldots ,p_d\}$ of primes. The space derivatives in the various
directions given by the various primes correspond to
Fermat quotient operators $\delta_{\cP}=\{\delta_{p_1},\ldots ,
\delta_{p_d}\}$ that are subjected to certain commutation
relations.

Let us denote by $\bZ_{(\cP)}$ the ring of fractions of $\bZ$ with respect
to the multiplicative system of all integers coprime to the primes in $\cP$.
Our first task is to define ``arithmetic PDEs''
with respect to $\delta_{\cP}$ on a smooth scheme $X$ over $\bZ_{(\cP)}$
of relative dimension $d_3\geq 1$, which will not necessarily be of a linear
nature. In the applications  $d_3$ will be taken to be $1$. We proceed
to explain the main ideas when in the scheme $X$ is affine.

Following the case of classical derivations \cite{kolchin, annals}, we
would like to define our PDEs as ``functions on jet spaces.'' Thus, we first
construct {\it arithmetic jet spaces} $\cJ^r_{\cP}(X)$ of order
$r=(r_1,\ldots ,r_d) \in \bZ_+^r$ (that is to say, of order $r_k$ with
respect to $\delta_{p_k}$ for all $k$) by analogy with the case of derivations.
When dealing with a single prime $p$, a $p$-adically completed version of
these spaces was introduced and thoroughly analyzed by the first author
in \cite{char, pjets}, and in a number of subsequent papers; cf. \cite{book}
for an exposition of this theory.
In the case where the number of primes $d$ is strictly greater than $1$,
our arithmetic jet spaces were also independently introduced in a recent
preprint by Borger \cite{borger}, where they were denoted by $W_{r*}(X)$.

As in turns out, in most cases the ring $A$ of global functions on a
jet space $\cJ^r_{\cP}(X)$
does not contain ``interesting elements'' (e.g. functions that qualify as
analogues of linear differential operators).
However, we show that, for each of the primes $p_k$ in $\cP$, we can
construct ``interesting'' formal functions  $f_{k} \in A\whk$, where $A\whk$
is the $p_k$-adic completion of $A$. Further, we would then like to be able
to say when certain families $f_{1},\ldots ,f_{d}$ of such formal
functions ``can be glued'' together. Since the ``domains of definition'' of
these $f_{k}$s are disjoint (each $f_{k}$ is defined on a
``vertical tubular neighborhood'' of  ${\rm Spec}\, A/p_kA$), we need
an ``indirect'' approach to the gluing of these formal functions.
We propose here to do this ``gluing'' by some sort of ``analytic
continuation'' between different primes, and
proceed as follows.

We first fix a $\bZ_{(\cP)}$-point $P$ of $X$. This point
can be lifted canonically to a $\bZ_{(\cP)}$-point $P^r$ of
$\cJ^r_{\cP}(X)$, and we still denote by $P^r \subset A$ the ideal of
the image of $P^r$. Then, we declare that a family
$f=(f_{k}) \in \prod_{k=1}^d A\whk$ can be {\it analytically continued
along $P$} if there exists $f_0 \in A^{\widehat{P^r}}$ such that,
for each $k=1,\ldots ,d$, the images of $f_0$ and $f_k$ in
$A^{\widehat{(p_k,P^r)}}$ coincide. Here, the superscripts ``$\widehat{\ }$''
mean completions with respect to the corresponding ideals.
Intuitively, $f_0$ can be viewed as a formal function defined on a
``horizontal tubular neighborhood'' of the canonical lift to the jet
space of our fixed point $P$ of $X$. This horizontal tubular neighborhood
meets each of the vertical tubular neighborhoods transversally in one point,
hence, in particular, the union of the vertical tubular neighborhoods
with the horizontal one is connected and, actually, it
``looks like a tree'' so is morally ``simply connected;'' this justifies
the terminology of {\it analytic continuation}.
Families $f$ that can be analytically continued along $P$ will be
referred to as {\it $\delta_{\cP}$-functions}, and
will be viewed as arithmetic (non-linear) PDEs on $X$.
The analytic continuation concept really depends on the choice of $P$. However,
we will provide a description of this dependence in the special cases
under consideration.

This can all be generalized to the case when $X$ is not affine, though some
extra difficulties have to be overcome. In the case when $X$ is a group
scheme $G$ and $P$ is the identity, the ``additive''
$\delta_{\cP}$-functions will be called
{\it $\delta_{\cP}$-characters}; morally these play the role of the
linear arithmetic PDEs on $G$.
Our main result is the determination
of all $\delta_{\cP}$-characters on $G=\bG_a$, $G=\bG_m$, or $G=E$.
(Here $E$ is, again, an elliptic curve.)
Indeed let $e_k=(0,\ldots ,1,\ldots ,0)\in \bZ_+^d$,
with $1$ on the $k$th place,  and let $e=e_1+\cdots +e_d=(1,\ldots ,1)$.
Then, in each of the $3$ cases $G=\bG_a$, $G=\bG_m$, or $G=E$, the spaces
of $\delta_{\cP}$-characters are generated by a fundamental
$\delta_{\cP}$-character of order $se=(s,\ldots ,s)$,
where $s=0,1,2$, respectively. These fundamental $\delta_{\cP}$-character
will be denoted by $\psi_{G}=\psi^0_a$, $\psi_{G}= \psi^s_m$, and
$\psi_{G}=\psi_E^{2e}$, respectively.
(As a general convention, upper indices for PDEs stand for order.)
In each of the $3$ cases we have decompositions
\begin{equation}
\label{appro}
\psi_G=(\bar{\psi}_{\bar{p}_k} \psi_{p_k})
\end{equation}
where $\psi_{p_k}$ come from $\delta_{p_k}$-characters
of $G$ of order $s=0,1,2$, respectively (in the sense of the arithmetic
ODE theory in \cite{char}) and $\bar{\psi}_{\bar{p}_k}$ come from
are $\delta_{\bar{p}_k}$-characters of $\bG_a$, where $\bar{p}_k$ is the
set $\cP \backslash \{p_k\}$.
Each $\bar{\psi}_{\bar{p}_k}$ is, loosely speaking, the product of
the ``symbols'' of all $\psi_{p_l}$ with $l \neq k$. These
``symbols'' are related to  the Euler factors of certain corresponding
$L$-functions. (The case $G=\bG_a$ does not actually require $L$-functions;
for $G=\bG_m$ and $G=E$ the $L$-functions are the Riemann zeta function
and the Hasse-Weil $L$-function respectively.) The equality (\ref{appro})
 will be viewed as a {\it Dirac decomposition} of $\psi_G$; cf.
the next subsection for analytic analogues of this.
Since one uses  \cite{char}, the construction of $\psi_G$
is global.
The fact that $\psi_G$ generates the space of $\delta_{\cP}$-characters
will be proved by a purely local method
based on results of Honda, Hill, and Hazewinkel; cf \cite{hazebook}.
We will also show, in each of the $3$ cases, that there are plenty of
$\delta_{\cP}$-characters that can be analytically continued along
{\it any} of the $\bZ_{(\cP)}$-points of $G$ (although this is not true
in general about $\psi_G$ itself). Finally we will also show, in each
of the $3$ cases,  that the solution group of $\psi_G=0$ in appropriate
cyclotomic rings consists of the torsion points of $G$ with values in
those rings; this should be viewed as an arithmetic analogue of an
analytic statement that says that the only
harmonic functions (in an appropriate  setting) are the constants; cf.
the discussion in the next subsection.

\subsection{Analytic analogues}
\label{anan}
When $d=2$ (and $G=\bG_m$ or $G=E$), our fundamental
$\delta_{\cP}$-characters $\psi_G$ are arithmetic analogues
of operators in analysis that are related to Laplacians. We explain here
this analogy.

Our fundamental $\delta_{\cP}$-character $\psi_m^{e}$ on $\bG_m$  should
be viewed as an arithmetic  analogue of
the partial differential operator
\begin{equation}
\label{trew}
\begin{array}{ccl}
\psi^{e}_{z\bar{z}}:C^{\infty}(D,\bC^{\times}) & \ra & C^{\infty}(D,\bC)\\
u & \mapsto & \psi^{e}_{z\bar{z}}(u):={\displaystyle \frac{1}{4} \Delta \log u=
\partial_z \partial_{\bar{z}} \log u},\end{array}
\end{equation}
  where $D \subset \bC$ is a domain, $z=x+iy$ is the complex
coordinate on  $D$,  and $\Delta=\partial_x^2+\partial_y^2$ is
the Euclidean Laplacian. (Here $\partial_x,\partial_y,\partial_z,
  \partial_{\bar{z}}$ are the corresponding partial derivative operators.)
  Note that, like our arithmetic $\psi_m^e$, the
operator $\psi_{z\bar{z}}^e$ is a group homomorphism
and has a  ``Dirac decomposition:''
$$\psi^{e}_{z\bar{z}}(u)=
 \partial_z \left(\frac{\partial_{\bar{z}} u}{u}\right)=
 \partial_{\bar{z}}\left(\frac{\partial_z u}
{u}\right)$$
which is analogous to the decomposition (\ref{appro}).
 Note that if we equip $D$ and $\bC^{\times}$ with their usual complex
structure, and we equip $\bC^{\times}$ with the conformal metric
  $\frac{dzd\bar{z}}{z\bar{z}}$  then $\psi^e_{z\bar{z}}$
    is an instance of the operator in \cite{jost}, p. 124 appearing in
the theory of harmonic maps between Riemann surfaces.
 Also recall from \cite{jost}, p. 63, that for real positive $u$ we have
 $\Delta \log u=-K_u \cdot u^2$,
 where $K_u$ is the curvature of the conformal metric $g=u^2dz d\bar{z}$.
If instead of  $D$ we take an arbitrary  Riemann surface $\Sigma$ then
$c_1(g):=-\frac{1}{4 \pi i} K_u \cdot u^2 dz \wedge d\bar{z}$ is the
Chern form of the metric $g$. So our character $\psi_m^{e}$  can also
be viewed as an arithmetic analogue of the operator
 \begin{equation}\label{pisu}
 \begin{array}{ccc}
\psi^{e}_{\Sigma}:\{\text{conformal metrics on $\Sigma$}\} & \ra &
\{(1,1)\text{-forms on $\Sigma$}\}\\
g & \mapsto & \psi^{e}_{\Sigma}(g)=c_1(g)
\end{array}\, .
 \end{equation}

Similarly, our fundamental $\delta_{\cP}$-character $\psi^{2e}_E$ on
an elliptic curve $E$  can be viewed as an arithmetic analogue of an
analytic operator that we now describe.
 Let $D$ be a domain contained in the upper half plane
$\{\tau=x+iy \in \bC\ ;\ y>0\}$. Let $E=(D \times \bC)/\sim$
where $(\tau_1,z_1) \sim (\tau_2,z_2)$ iff $\tau_1=\tau_2=:\tau$
and  $z_2-z_1 \in \bZ\tau+\bZ$. Let $\pi:E \ra D$ be the
canonical projection
(so $E$ is the universal elliptic curve over $D$).
Let $C^{\infty}_{\pi}(D,E)$ the set of all maps $\sigma\in C^{\infty}(D,E)$
such that $\pi \circ \sigma=1_D$. Let $\log_E:E \ra \bC$ be
the (multivalued) map obtained by composing
the (multivalued) inverse of the canonical projection
$D \times \bC \ra E$ with the second projection
$D \times \bC \ra \bC$. Then  our arithmetic $\psi^{2e}_E$ can be
 viewed as an analogue of the order 4 map
\begin{equation}
\label{trewqq}
\begin{array}{ccl}
\psi^{2e}_{\tau\bar{\tau}}:C^{\infty}_{\pi}(D,E) & \ra & C^{\infty}(D,\bC)\\
u & \mapsto & \psi^{2e}_{\tau\bar{\tau}}(u):={\displaystyle
\frac{1}{16} \Delta \Delta \log_E u=\partial_{\tau}^2
\partial^2_{\bar{\tau}}  \log_E u},\end{array}
\end{equation}
where $\Delta=\partial_x^2+\partial^2_y$ is the Euclidean Laplacian.
(In order to see that this map is well defined, notice that if
$\varphi:D \ra \bC$ is a smooth map satisfying
 $\varphi(\tau)\in \bZ \tau +\bZ$ for all $\tau \in D$, then there
exist $m,n \in \bZ$ such that  $\varphi(\tau)=m\tau+n$ for all $\tau \in D$
and hence $\varphi$ is in the kernel of
$\partial_{\tau}^2 \partial^2_{\bar{\tau}}$.)
The map $\psi^{2e}_{\tau\bar{\tau}}$ is a group homomorphism, and has
a ``Dirac decomposition''
\begin{equation}
\label{ddiirr}
\psi^{2e}_{\tau\bar{\tau}}(u)=\partial_{\bar{\tau}}^2 (\mu_{\tau}(u))=
\partial_{\tau}^2 (\mu_{\bar{\tau}}(u))\, ,
\end{equation}
where
\begin{equation}
\label{maninan}
\mu_{\tau}(u)=\partial_{\tau}^2  \log_E u\, ,\quad
\mu_{\bar{\tau}}(u)=\partial_{\bar{\tau}}^2  \log_E u\, .
\end{equation}
This is again analogous to (\ref{appro}).
The map $\mu_{\tau}$ is the analytic expression of  the Manin map
in \cite{manin}, and is well defined because any
 map $\varphi(\tau)=m\tau+n$ ($m,n \in \bZ$)
 is in the kernel of $\partial_{\tau}^2$. Also,
 $\mu_{\bar{\tau}}$ is well defined because any $\varphi$ as above
is in the kernel of $\partial_{\bar{\tau}}^2$.
Note also that one could  consider, in this analytic setting, the
order $2$ map
\begin{equation}
\begin{array}{ccl}
\psi^e_{\tau\bar{\tau}}:C^{\infty}_{\pi}(D,E) & \ra & C^{\infty}(D,\bC)\\
u & \mapsto & \psi^e_{\tau \bar{\tau}}(u):={\displaystyle
\frac{1}{4}\Delta \log_E u=\partial_{\tau}\partial_{\bar{\tau}}
\log_E u}\, .
\end{array}
\end{equation}
(This map is well defined because any  $\varphi(\tau)=m\tau+n$ is in the
 kernel of $\partial_{\tau} \partial_{\bar{\tau}}$.) The map
$\psi^2_{\tau\bar{\tau}}$ does not seem to have, however, a natural
``Dirac decomposition''
and it has no arithmetic analogue. (The ``obvious" candidate for a
``Dirac decomposition'' for $\psi^e_{\tau\bar{\tau}}$, analogous
 to (\ref{ddiirr}), fails to make sense because
the operator $u \mapsto \partial_{\tau} \log_E u$ is not well defined.)

 Let us recall that the  Manin map \cite{manin} is an order $2$ ODE map
introduced by Manin  for abelian varieties, in order to prove the Mordell
conjecture over function fields. A different construction of this map
was given in \cite{annals}.
The construction in \cite{annals} has an arithmetic analogue which
was introduced in \cite{char}, giving rise to an order $2$ ODE
{\it arithmetic Manin map} $\psi^2_{p_k}$. Our  arithmetic PDE
$\psi^{2e}_E$ has a Dirac decomposition involving the operators $\psi_{p_k}^2$.

Instead of the family $E \ra D$ above we could as well consider a fixed
elliptic curve $E_0=\bC/(\bZ\tau_0+\bZ)$. If $\log_{E_0}:E_0 \ra \bC$ is
the (multivalued) inverse of the canonical projection $\bC \ra E_0$,
 and $D$ is a domain in the complex plane (with coordinate $w=x+iy$)
then we have a well defined map
\begin{equation}
\begin{array}{ccl}
\psi^e_{w\bar{w}}:C^{\infty}(D,E_0) & \ra & C^{\infty}(D,\bC)\\
u & \mapsto & \psi^e_{w\bar{w}}(u):={\displaystyle
\frac{1}{4}\Delta \log_{E_0} u=\partial_{w}\partial_{\bar{w}}
\log_{E_0} u}\, ,
\end{array}
\end{equation}
admitting an obvious ``Dirac decomposition.'' There is an arithmetic
analogue of this map
which we are not going to discuss in the present paper because such a
discussion would require developing our theory in a setting slightly more
 general than the one we have adopted here. In this more general setting,
$\bQ$ needs to be replaced by a number field $F$ containing an imaginary
quadratic  field $K$. Then the analogue of $E_0$ is an elliptic curve over
$F$ with complex multiplication by $K$, the primes in $\cP$ are replaced by
primes of $F$ of degree one, where $E_0$ has ordinary good reduction, and
the $L$-function appearing is the $L$-function associated to the
corresponding Grossencharacter.

\subsection{$p$-adic notation}
\label{cou}
Given a prime $p\in \bZ$, we set
$$\begin{array}{ccl}
\bZ_{(p)} & = & \text{localization of $\bZ$ at $(p)$}\, ,\\
\bZ_p & = &   \text{ring of $p$-adic integers}\, ,
\end{array}$$
notation that we use throughout the  paper.
Given a family $\cP=\{p_1, \ldots ,p_d\}$ of primes in $\bZ$, we
set
$$\begin{array}{ccl}
\bZ_{(\cP)} & = & \bigcap_{k=1}^d \bZ_{(p_k)}.
\end{array}$$
For each $k=1,\ldots ,d$, we set
$$
\bar{p}_k=\cP \backslash \{p_k\}=\{p_1,\ldots ,p_{k-1},p_{k+1},
\ldots,p_d\}\, .
$$

If $A$ is a ring and $I$ is an ideal, we let
$A^{\widehat{I}}$ denote the completion of $A$ with respect to  $I$.
If the ideal $I$ is generated by a family of elements $a$, we shall
write $A^{\widehat{a}}$ instead of $A^{\widehat{I}}$.
In particular, for any prime $p\in \bZ$, we have that
$A^{\widehat{p}}$ is the $p$-adic completion of $A$.

If $X$ is any Noetherian scheme, we
denote by $X^{\widehat{p}}$ the $p$-adic completion of $X$. For
any ring $A$ and family of primes $\cP=\{p_1, \ldots ,p_d\}$ as above,
we set
$$\begin{array}{rcl}
A_{\cP}& = & \prod_{k=1}^d A\whk.\end{array}
$$

\subsection{Structure of the paper}
In \S 2, we present our main concepts, in particular, the notions of
$\delta_{\cP}$-jet spaces and $\delta_{\cP}$-characters. In \S 3, we
state and prove our main results determining all the
$\delta_{\cP}$-characters on one dimensional groups.
In \S 4 we present a number of remarks and open questions.
\medskip

\subsection*{Acknowledgement}
The authors are indebted to J. Borger for discussions on Witt vectors
and $\lambda$-rings, and to D. Bertrand for correspondence on transcendence.

\section{Main concepts}
\setcounter{theorem}{0}

\subsection{$\12$-rings}
For any prime $p \in \bZ$, we consider the polynomial
$$C_p(X,Y):=\frac{X^p+Y^p-(X+Y)^p}{p} \in \bZ[X,Y]\, .$$
Let $A$ be a ring and $B$ an $A$-algebra.
For $a \in A$, we denote the element $a \cdot 1_B \in B$ by $a$ also.

\begin{definition}
\label{pder}
A map $\delta_p:A \ra B$ is a {\it $p$-derivation} if
\begin{equation}
\begin{array}{rcl}
\delta_p(a+b) & = & \delta_p a+\delta_p b + C_p(a,b)\, ,\\
\delta_p(ab) & = & a^p \delta_p b+b^p \delta_p a + p \delta_p a
\delta_p b \, ,
\end{array}
\end{equation}
for all $a,b \in A$. \qed
\end{definition}

If $\delta_p$ is a $p$-derivation, then the map
\begin{equation}
\label{phi}
\begin{array}{rcl}
A & \stackrel{\phi_p}{\ra} & B \\
a & \mapsto & \phi_p(a):= a^p+p\delta_p a
\end{array}
\end{equation}
is a ring homomorphism, and we have that
\begin{equation}
\label{cong}
\phi_p(a) \equiv a^p\quad {\rm mod}\; p
\end{equation}
for all $a \in A$.

Conversely, let $\phi_p:A \ra B$ be a ring homomorphism
satisfying (\ref{cong}), which we refer to as a {\it lift of Frobenius}
homomorphism. Let us assume in addition that  $p$ is a non-zero divisor in $B$.
Then $\phi_p$ has the form (\ref{phi}) for a unique $p$-derivation $\delta_p$,
$$
\delta_p a=\frac{\phi(a)-a^p}{p}\, ,
$$
and we have that
$$
\phi_p\delta_pa=\delta_p\phi_p a
$$
for all $a \in A$. We say that $\delta_p$ and $\phi_p$ are {\it associated}
to each other.

In particular, the ring $A=\bZ$ has a unique $p$-derivation
$\delta_p:\bZ  \ra \bZ$, given by
$$
\delta_p a := {\displaystyle \frac{a-a^p}{p}}\, .
$$

If $p_1$ and $p_2$ are two distinct primes in $\bZ$, we now consider the
polynomial $C_{p_1,p_2}$ in the ring $\bZ[X_0,X_1,X_2]$ defined by
\begin{equation}
\label{commutator}
C_{p_1,p_2}(X_0,X_1,X_2) := \frac{C_{p_2}(X_0^{p_1},p_1X_1)}{p_1}
-\frac{C_{p_1}(X_0^{p_2},p_2 X_2)}{p_2}
 -\frac{\delta_{p_1} p_2}{p_2} X_2^{p_1}+ \frac{\delta_{p_2}p_1}{p_1}
X_1^{p_2}\, .
\end{equation}
Let us notice that
\begin{equation}
\label{square}
C_{p_1,p_2} \in (X_0,X_1,X_2)^{\min\{p_1,p_2\}}\, .
\end{equation}

\begin{definition}
\label{p1p2}
Let $\cP=\{p_1,\ldots ,p_d\}$ be a finite set of  primes in $\bZ$.
A $\12$-{\it ring} is a ring $A$  equipped with
$p_k$-derivations $\delta_{p_k}:A \ra A$,  $k=1,\ldots ,d$, such that
\begin{equation}
\label{identity}
\delta_{p_k}\delta_{p_l}a-\delta_{p_l}\delta_{p_k}a=C_{p_k,p_l}(a,
\delta_{p_k}a,\delta_{p_l}a)
\end{equation}
for all $a \in A$, $k,l=1,\ldots ,d$.
A {\it homomorphism
of $\delta_{\cP}$-rings} $A$ and $B$ is a homomorphism
of rings $\varphi: A \rightarrow B$ that commutes
with the $p_k$-derivations in $A$ and $B$, respectively.
\qed
\end{definition}

\medskip

If $\phi_{p_k}$ is the homomorphism (\ref{phi}) associated to $\delta_{p_k}$,
condition
(\ref{identity}) implies that
\begin{equation}
\label{commu}
\phi_{p_k}\phi_{p_l}(a)=\phi_{p_l}\phi_{p_k}(a)
\end{equation}
for all $a \in A$. Conversely, if the commutation relations (\ref{commu})
hold, and the $p_k$s are non-zero divisors in $A$, then conditions
(\ref{identity}) hold, and we have that
$$
\phi_{p_k}\delta_{p_l} a  =  \delta_{p_l}\phi_{p_k}a
$$
for all $a \in A$.

\medskip

\begin{remark}
The concept of a ``lift of Frobenius'' homomorphism that is the basis for
the definitions given above, is classical and goes  back
to work of Frobenius, Chebotarev, and Artin on number fields.
It plays a key role in the theory of Witt vectors (in particular, in
crystalline cohomology), and it resurfaced in the context of $K$-theory
through Grothendieck's concept of lambda ring;
cf. \cite{groth, hazebook, wilk, joy, bor} for  this circle of ideas.
This concept was taken in \cite{char, difmod, book} as the starting point
for developing an arithmetic analogue of ordinary differential equations
on commutative algebraic groups and on moduli spaces (such as modular curves,
Shimura curves, and Siegel modular varieties). In this article,
we aim at extending these ideas to the partial differential case
(at least for one dimensional algebraic groups).
\qed
\end{remark}

\medskip

\begin{remark}
\label{comaredel}
Let $A$ be a $\12$-ring. Then for all $k$, the $p_k$-adic completions
$A\whk$ are $\12$-rings in a natural way. For $\phi_{p_l}$
extends to $A^{\widehat{p_k}}$ by continuity. The
condition $\phi_{p_k}a\equiv a^{p_k}$ mod $p_k$ in $A\whk$ holds
by continuity, while
the condition $\phi_{p_l}a\equiv a^{p_l}$ mod $p_l$ in $A\whk$ holds
because $p_l$ is invertible in $A^{\widehat{p_k}}$.
\qed
\end{remark}

\medskip

We let $\bZ_+=\{0,1,2,3, \ldots \}$, and
let $\bZ_+^d$ be given the product order, defined by
declaring that
$$
r=(r_1,\ldots ,r_d) \leq r'=(r'_1,\ldots ,r'_d)
$$
if $r_k\leq r'_k$ for all $k=1,\ldots ,d$. The weight $|r|$ of the
$d$-tuple $r$ is defined by $|r|=\sum_{j=1}^d r_j$. Finally, we
let $e_k$ be the element of $\bZ_{+}^d$ all of whose components are
zero except the $k$-th, which is $1$.

In more generality, we have the following:

\begin{definition}
\label{prol}
A {\it $\delta_{\cP}$-prolongation system} $A^*=(A^r)$ is an inductive
system of rings $A^r$ indexed by $r \in \bZ_+^d$, provided with
transition maps $\varphi_{rr'}:A^r \ra A^{r'}$ for any pair of indices
$r$, $r^{'}$ such that $r \leq r'$, and equipped with $p_k$-derivations
$$
\delta_{p_k}:A^r\ra A^{r+e_k}\, ,
$$
$k=1,\ldots ,d$, such that (\ref{identity}) holds for all $k$, $l$, and
such that
$$
\varphi_{r+e_k,r'+e_k}\circ \delta_{p_k}=\delta_{p_k} \circ
\varphi_{rr'} :A^r \ra A^{r'+e_k}
$$
for  all $r\leq r'$ and all $k$.
A morphism of prolongation systems $A^* \ra B^*$ is a system of ring
homomorphisms
$u^r:A^r \ra B^r$ that commute with the
$\varphi$s and the $\delta$s of $A^*$ and $B^*$, respectively.
\qed
\end{definition}
\medskip

Any $\delta_{\cP}$-ring $A$ induces a $\delta_{\cP}$-prolongation system
$A^*$ where $A^r=A$ for all $r$ and $\varphi=$identity. If $A$ is a
$\delta_{\cP}$-ring and $A^*$ is the associated
$\delta_{\cP}$-prolongation system, we say that
a $\delta_{\cP}$-prolongation system $B^*$ is a
{\it $\delta_{\cP}$-prolongation system over $A$}
if it is equipped with a morphism $A^* \ra B^*$.
We have a natural concept of morphism of $\delta_{\cP}$-prolongation
systems over $A$.

\begin{example}
\label{basicex}
Let  $S \subset \bZ$ be a multiplicative system of integers coprime
to $p_1, \ldots, p_d$, and let $\bZ_S=S^{-1}\bZ$ be the ring of fractions.
Given an integer $m$ coprime to $p_1, \ldots, p_d$, let
$\zeta_m=\exp \left(\frac{2 \pi \sqrt{-1}}{m}\right)$.
We consider the ring $A=\bZ_S[\zeta_m]$.
Let $\phi_{p_1},\ldots ,\phi_{p_d} \in G(\bQ(\zeta_m)/\bQ)=
(\bZ/m\bZ)^{\times}$ be the Galois elements corresponding
to the classes of $p_1,\ldots ,p_d$, respectively. Then
$\phi_{p_k}(a) \equiv a^{p_k}$ mod $p_k$ for $k=1,\ldots ,d$ and $a \in A$.
Thus, $A$ is a $\12$-ring with respect to the $p_k$-derivations
$\delta_{p_k}$ associated to the $\phi_{p_k}$s.
\qed
\end{example}
\medskip

{\it From now on,
$A_0$ will denote  a $\delta_{\cP}$-ring that is an integral Noetherian domain
of characteristic zero, with fraction field $K_0$.
Notice that  $K_0$ has a naturally induced structure
of $\delta_{\cP}$-ring.}
\medskip

\begin{example}
Let $A=A_0[[q]]$ be the power series ring in the indeterminate $q$.
We let
$$
\phi_{p_k} : A \ra A
$$
$1\leq k \leq d$ be the family of homomorphisms defined by
$$\phi_{p_k}(\sum c_n q^n)=\sum \phi_{p_k}(c_n) q^{np_k}\, .$$
$k=1, \ldots, d$. Then, as before
 $A$ is a $\12$-ring with respect to the $p_k$-derivations
  $\delta_{p_k}$ associated to the $\phi_{p_k}$s.
\end{example}

\medskip

Let $x$ be an $n$-tuple of variables.
We consider $n$-tuples of variables $x_{i}$ indexed by vectors
$i=(i_1,\ldots ,i_d)$ in $\bZ_+^d$.
 We set
$x=x_{(0,\ldots, 0)}$, $\delta_{\cP}^i=\delta_{p_1}^{i_1}\ldots
\delta_{p_d}^{i_d}$, and $\phi_{\cP}^i=\phi_{p_1}^{i_1}\ldots
\phi_{p_d}^{i_d}$.

Let $K_0\{x\}$ be the ring of polynomials
$$
K_0\{x\}:=K_0[x_{i} :\, i \in \bZ_+^d]
$$
with $K_0$-coefficients in the variables
$x_{i}$. We extend the homomorphisms $\phi_{p_k}:A_0 \ra A_0$ to ring
endomorphisms $\phi_{p_k}:K_0\{x\}\ra K_0\{x\}$ by the formulae
$$
\phi_{p_k}(x_{i})  =  x_{i+e_k}\, ,
$$
 so that
$x_i=\phi_{\cP}^i x$ for all $i$.
Clearly $\phi_{p_k}\phi_{p_l}(a)=\phi_{p_l}\phi_{p_k}(a)$
for all  $a \in K_0\{x\}$, and all $k,l$. If we consider the  $p_k$-derivations
$
\delta_{p_k}:K_0\{x\} \ra K_0\{x\}
$
associated to the $\phi_{p_k}$s
then $K_0\{x\}$ is a $\12$-ring. Notice that $K_0\{x\}$ is generated
as a $K_0$-algebra by the elements $\delta_{\cP}^i x$, $i \in \bZ_+^d$:
$$
K_0\{x\}=K_0[\delta_{\cP}^i x : \, i \in \bZ_+^d]\, .
$$

\begin{example}
\label{moo}
We define the ring of $\12$-{\it polynomials} $A_0\{x\}$ to be the
$A_0$-subalgebra of $K_0\{x\}$ generated by all the elements
$\delta_{\cP}^i x$:
$$
A_0\{x\}:=A_0[\delta_{\cP}^i x :\, i \in \bZ_+^d]\, .
$$
Notice that $A_0\{x\}$ is strictly larger than the ring
$A_0[x_{i} :\, i\in \bZ_+^d]$. And notice also that the family
$\{\delta_{\cP}^i :\, i \in \bZ_+^d\}$ is algebraically independent
over $A_0$, so $A_0\{x\}$ is a ring of polynomials in the
``variables'' $\delta_{\cP}^i x$.

The ring $A_0\{x\}$ has a natural structure of  $\12$-ring due to the
following:

\begin{lemma} We have that
$\delta_{p_k} A_0\{x\} \subset A_0\{x\}$ for $k=1,\ldots ,d$.
\end{lemma}

{\it Proof}. By (\ref{pder}), the sets
$S_k:=\{a \in K_0\{x\} :\, \delta_{p_k}a \in A_0\{x\}\}$
are $A_0$-subalgebras of $K_0\{x\}$, so it is enough to show that
$\delta_{\cP}^i x \in S_k$ for all $i$ and $k$.
We use the commutation relations ($\ref{identity}$) to check this by
induction on $(i,k) \in \bZ_+^d\times \bZ_+$ with respect to the
lexicographic order.
\qed
\end{example}

\begin{example}
\label{goo}
Proceeding exactly as in Example \ref{moo}, the system
$$
A_0[\delta_{\cP}^i x :\, i \leq r]
$$
has a natural structure of $\delta_{\cP}$-prolongation system. \qed
\end{example}

\begin{example}
Let $T$ be a tuple of indeterminates, and let
$$
A^r=A_0[[\delta_{\cP}^i T :\, i\leq r]]\, .
$$
Then, the structure of $\12$-prolongation sequence in Example \ref{goo}
induces a structure of $\12$-prolongation sequence on the sequence of rings
$(A^r)$. Indeed, the $p_k$-derivation
$\delta_{p_k}$ sends the ideal
$$
I_r:=(\delta_{\cP}^i T : \, i \leq r) \subset A^r
$$
into the ideal $I_{r+e_k}\subset A^{r+e_k}$; cf. (\ref{square}). \qed
\end{example}

\subsection{$\delta_{\cP}$-jet spaces}
As in the case of a single prime \cite{char}, we now have the following
existence result for a universal prolongation sequence.

\begin{proposition}
\label{UP}
Let $A^0$ be a finitely generated $A_0$-algebra. Then there exists a
$\delta_{\cP}$-prolongation sequence $A^*$ over $A_0$, with $A^r$
finitely generated over $A_0$,
satisfying  the following property: for any $\delta_{\cP}$-prolongation
sysytem $B^*$ over $A_0$ and any $A_0$-algebra homomorphism
$u:A^0 \ra B^0$, there exists a unique morphism of
$\delta_{\cP}$-prolongation systems $u^*:A^* \ra B^*$
such that $u^0=u$.
\end{proposition}

{\it Proof}. We express the finitely generated algebra $A^0$ as
$$
A^0=\frac{A_0[x]}{(f)}
$$
for a tuple of indeterminates $x$, and a tuple of polynomials $f$.
Then we set
$$
A^r=\frac{A_0[\delta_{\cP}x^i : \, i\leq r]}{(\delta_{\cP}^i f :\,
i\leq r)}\, .
$$
Using  Example \ref{goo} and the identities (\ref{identity}), we check easily
that $A^*=(A^r)$ has the universality property in the statement. \qed

\begin{definition}
\label{arde}
Let $X$ be an affine scheme of finite type over $A_0$. Let $A^0=\cO(X)$
and let $A^r$ be as in Proposition \ref{UP}. Then the scheme
$${\mathcal J}_{\cP}^r(X):={\rm Spec}\, A^r$$
is called the {\it $\delta_{\cP}$-jet space of order $r$} of $X$. \qed
\end{definition}

By the universality property in Proposition \ref{UP},
up to isomorphism the scheme ${\mathcal J}_{\cP}^r(X)$ depends
on $X$ alone, and is functorial in $X$: for any morphism $X \ra Y$ of
affine schemes of finite type, there are induced morphisms of schemes
$$
{\mathcal J}_{\cP}^r(Y) \ra  {\mathcal J}_{\cP}^r(X)\, .
$$

\begin{remark}
In the case when $\cP$ consists of a single prime $p$, the $p$-adic
completions of the schemes $\cJ^r_p(X)$ were introduced and thoroughly
studied in \cite{char,pjets,book}. For arbitrary $\cP$,
the schemes $\cJ^r_{\cP}(X)$ above were independently introduced by Borger
in \cite{borger}, where they are denoted by $W_{r*}(X)$.
\end{remark}

\begin{lemma}
\label{dorinpiept}
Let $X$ be an affine scheme of finite type over $A_0$ and let $Y \subset X$ be
a principal open set of $X$, $\cO(Y)=\cO(X)_f$. Then
$\cO({\mathcal J}_{\cP}^r(Y)) \simeq \cO({\mathcal J}_{\cP}^r(X))_{f_r}$
where $f_r=\prod_{i\leq r}\phi_{\cP}^i(f)$.
In particular, the induced morphism ${\mathcal J}_{\cP}^r(Y) \ra
{\mathcal J}_{\cP}^r(X)$ is an open immersion
whose image is principal, and if we view this morphism
as an inclusion and $Z \subset X$ is another principal open set, then we have
that
$$
{\mathcal J}_{\cP}^r(Y\cap Z)={\mathcal J}_{\cP}^r(Y)\cap
{\mathcal J}_{\cP}^r(Z)\, .
$$
\end{lemma}

{\it Proof}. We can check easily that $\cO({\mathcal J}_{\cP}^r(X))_{f_r}$
has the universality property of $\cO({\mathcal J}_{\cP}^r(Y))$.
The $\delta_{p_k}$-derivations on $\cO({\mathcal J}_{\cP}^r(X))_{f_r}$
are defined via the formula
$$
\delta_{p_k}\left(\frac{a}{b}\right) =\frac{b^{p_k}\delta_{p_k} a-a^{p_k}
\delta_{p_k}b}{b^p\phi_{p_k}(b)}\, .
$$
\qed

\begin{definition}
\label{lkj}
Let $X$ be a scheme of finite type over $A_0$. An affine open covering
\begin{equation}
\label{dec}
X=\bigcup_{i=1}^m X_i
\end{equation}
is called {\it principal} if
\begin{equation}
\label{princ}
\text{$X_i \cap X_j$ is principal in both, $X_i$ and $X_j$,}
\end{equation}
for all $i,j=1, \ldots ,m$.

Let $c=\{ X_i\}_{i=1}^m$ be a principal
covering of $X$. We define the {\it $\delta_{\cP}$-jet space of order $r$
of $X$ with respect to the covering $c$} by gluing ${\mathcal J}_{\cP}^r(X_i)$
along ${\mathcal J}_{\cP}^r(X_i\cap X_j)$, and
denote this jet space by ${\mathcal J}_{c,\cP}^r(X)$. \qed
\end{definition}
\medskip

In the case when $X$ is affine, we may use the
{\it tautological covering} $c$ of $X$, that is to say, the covering of $X$
with a single open set, the set $X$ itself. Then the space
$\cJ^r_{c,\cP}(X)$ coincides with the space $\cJ^r_{\cP}(X)$ of
Definition \ref{arde}.

Notice that any quasi-projective scheme $X$ admits a principal covering.

\begin{remark}
Generally speaking, the scheme ${\mathcal J}_{c,\cP}^r(X)$ in
Definition \ref{lkj} depends on the covering $c$. For instance,
if $X={\rm Spec}\, A_0[x]$ is the affine line, we may consider the principal
cover $\tilde{c}$
consisting of the two open sets $X_0={\rm Spec}\, A_0[x,1/x]$ and
$X_1={\rm Spec}\, A_0[x,1/(x-1)]$, respectively. Then, if
$\cP=\{p\}$, the scheme ${\mathcal J}_{\tilde{c},\cP}^1(X)$
is the union
\begin{equation}
\label{union}
{\rm Spec}\,  A_0\left[x,\delta_p x,\frac{1}{x (x^p+p\delta_p x)}\right]
\cup {\rm Spec}\,  A_0\left[x,\delta_p x,\frac{1}{(x-1)
(x^p+p\delta_p x-1)}\right]
\end{equation}
inside the plane
\begin{equation}
\label{plane}
{\rm Spec}\, A_0[x,\delta_p x]\, ,
\end{equation}
while the scheme ${\mathcal J}_{\cP}^1(X)$ corresponding to
the tautological covering $c$ of $X$ is the whole plane (\ref{plane}).
But the union (\ref{union}) does not equal this plane. For the
$K_0$-point of (\ref{plane}) with coordinates
$x=0$, $\delta_p x=1/p$ does not belong to (\ref{union}).

This phenomenon is similar to problems encountered in
\cite{CH,hrushovski,borger} and, as in loc.cit., there are
ways to fix this anomaly at the cost of developing more technology.
For our purposes here, these covering dependent rudimentary definition of
jet spaces will be sufficient.
Indeed some of the basic rings of functions that we are going to consider
will be independent
of the coverings; cf. Proposition \ref{GO}.
However, let us observe that the endofunctor $X \mapsto \cJ^r_{\cP}(X)$ on
the category of affine schemes of finite type over $A_0$
that we consider here has been extended in \cite{borger} to an endofunctor
on the category of algebraic spaces over $A_0$.
\qed
\end{remark}

\begin{remark}
If $X$ is a closed subscheme of a projective space over $A_0$, then $X$ has
a natural principal covering $c$ where the open sets are the complements of
the coordinate hyperplanes. We may thus attach to $X$ the scheme
${\mathcal J}_{c,\cP}^r(X)$ corresponding to this cover. In general, such
a definition depends upon the embedding of $X$ into a projective space.
\qed
\end{remark}

\begin{remark}
When $\cP=\{p\}$ consists of a single prime $p$, then the $p$-adic
completion ${\mathcal J}^r_{c,\cP}(X)^{\widehat{p}}$
(cf. the notation in \S \ref{cou}) coincides with the $p$-jet
space $J^r_p(X)$ over $\bZ_p$ in \cite{char,book}, and therefore, it
depends on $X$ only and not on the covering $c$. For
$\cP=\{p_1,\ldots ,p_d\}$, we obtain natural morphisms
$$
{\mathcal J}^{je_k}_{c,\cP}(X)\whk \simeq  \cJ^j_{c,\{p_k\}}(X)
$$
for each $j\in \bZ_{+}$, and consequently, natural
morphisms
\begin{equation}
\label{pisacraz}
{\mathcal J}^r_{c,\cP}(X)\whk \ra
\cJ^{r_ke_k}_{c,\cP}(X)\whk=
{\mathcal J}^{r_k}_{c,\{p_k\}}(X)\whk=J^{r_k}_{p_k}(X)\, .
\end{equation}
\qed
\end{remark}

\begin{remark}
\label{nouup}
Let $X$ be an affine scheme with its tautological cover $c$. Then the
system $(\cO(\cJ^r(X))\whk)$ has a natural structure of
$\delta_{\cP}$-prolongation system. It has a universality property
analogue to that in Proposition \ref{UP} with respect to
$\delta_{\cP}$-prolongation systems $B^*$ of $p_k$-adically complete
rings. \qed
\end{remark}

We have the following structure result for $p_k$-adic completions
of $\delta_{\cP}$-jet spaces.

\begin{lemma}
\label{YO}
Let us assume that $X={\rm Spec}\, A^0$ is affine over $A_0$, and let
$A_0[T]\ra A^0$ be an \'{e}tale map, where $T$ is a tuple of variables.
Consider $X$ with its tautological cover. Then there is a natural isomorphism
\begin{equation}
\label{oega}
\cO(\cJ^r_{\cP}(X))\whk \simeq (\otimes_{i\leq r-r_ke_k} A_i^0)
[\delta_{p_k}^j\phi_{\bar{p}_k}^iT : \, j\geq 1,\; i+je_k\leq r]\whk\, ,
\end{equation}
where $A_i^0 = A^0$, and $\otimes=\otimes_{A_0}$.
\end{lemma}

{\it Proof}. If $A^0=A_0[x]/(f)$, we set
$A_i^0:=A_0[\phi^i_{\cP} x]/(\phi^i_{\cP} f)$. Then we can check that
the right hand side of (\ref{oega}) has the universality property of the
left hand side, as explained in Remark \ref{nouup}. For this, we need to use
the fact that $p_l$ is invertible in any $p_k$-adic completion, $l \neq k$,
and use the argument in Proposition 3.13 in \cite{book}.
\qed

\begin{proposition}
\label{GO}
Let us assume that $X$ is a smooth scheme of finite type, quasi-projective
and with connected geometric fibers over $A_0=\bZ_{(\cP)}$. Let $c$ and $c'$ be
two principal coverings of $X$, and $\cJ^r_{c,\cP}(X)$ and $\cJ^r_{c',\cP}(X)$
be the corresponding jet spaces. Then, there is a natural isomorphism
$$
\cO(\cJ^r_{c,\cP}(X)\whk) \simeq \cO(\cJ^r_{c',\cP}(X)\whk)\, .
$$
\end{proposition}
\medskip

Thus, we may drop the covering from the notation and denote the
isomorphism class of rings $\cO(\cJ^r_{c,\cP}(X)\whk)$ simply
by $\cO(\cJ^r_{\cP}(X)\whk)$.

{\it Proof}. Let us assume that $X={\rm Spec}\, A^0$ has an \'{e}tale map
to an affine space over $A_0$, and consider a
covering $X=\cup X_j$ where $X_j={\rm Spec}\, A^0_{f_j}$. We let
$Y=\cJ^r_{\cP}(X)$ and $Y_j=\cJ^r_{\cP}(X_j)$
be the jet spaces corresponding to the tautological coverings
of $X$ and $X_j$, respectively, and set $U=\cup Y_j$.
We claim that the restriction map $\cO(Y\whk) \ra \cO(U\whk)$ is an
isomorphism, fact that is easily seen to imply the statement of the
Proposition.

For the proof of this claim, it is enough to check that the map
$\cO(\overline{Y}) \ra \cO(\overline{U})$ is an isomorphism,
where $\overline{Y}:=Y \otimes \bF_{p_k}$ and
$\overline{U}:=U \otimes \bF_{p_k}$, respectively.
Observe that $\overline{Y}$ is smooth over $\bF_p$ (cf. Lemma \ref{YO}),
and that $\overline{U}$ is an open subset of $\overline{Y}$
(cf. Lemma \ref{dorinpiept}). Thus, it suffices to check that
$\overline{Y}\backslash\overline{U}$ has codimension $\geq 2$ in
$\overline{Y}$. By Lemmas \ref{YO} and \ref{dorinpiept},
 we have identifications $\overline{Y}=\overline{X}^a \times \bA^b$ and
$\overline{Y_j}=\overline{X_j}^a \times \bA^b$, where $\bA^b$ is the
affine space of dimension  $b$ over $\bF_{p_k}$, and $\overline{X}^a$
and $\overline{X}^a_j$ are the $a$-fold products of $\overline{X}$ and $X_j$,
respectively. Let us set $H_j=\overline{X}\backslash \overline{X_j}$.
We may assume that all the $H_j$s are non-empty and different from
$\overline{X}$. So each $H_j$ is a hypersurface in $\overline{X}$. We have
that
$$
\overline{X}^a \backslash (\cup \overline{X}_j^a)\subset
\bigcup_{1 \leq i_1,\ldots ,i_a \leq a} S_{i_1,\ldots ,i_a}$$
where
$$
S_{i_1,\ldots, i_a}=\{(P_1,\ldots ,P_a)\in \overline{X}^a\, :
\;  P_{i_k} \in H_k,\; 1\leq k \leq a\}\, .
$$
Now, if $i_1=\cdots =i_a$ then $S_{i_1,\ldots ,i_a}=\emptyset$. On the
other hand, if at least two of the indices $i_1,\ldots ,i_a$ are different,
then $S_{i_1,\ldots ,i_a}$
has codimension $\geq 2$ in $\overline{X}^a$. This implies that
$\overline{Y} \backslash  \overline{U}$ has codimension $\geq 2$ in
$\overline{Y}$, as desired. \qed
\medskip

\subsection{$\delta_{\cP}$-functions}
We would like to see now how we can ``glue'' elements in various completions
of a given ring. We achieve this through the following definitions.
\medskip

\begin{definition}
\label{ideeea}
Let $A$ be a ring, $I$ be an ideal in $A$, and $\cP=\{p_1,\ldots ,p_d\}$
be a finite set of primes in $\bZ$ that are non-invertible in $A/I$.
We say that a family
\begin{equation}
\label{fma}
f=(f_k) \in \prod_{k=1}^d A^{\widehat{p_k}}
\end{equation}
can be {\it analytically continued along} $I$ if there exists
$f_0 \in A^{\widehat{I}}$ such that
the images of $f_0$ and $f_k$ in the ring
$A^{\widehat{(p_k,I)}}$ coincide for each $k=1,\ldots ,d$.
In that case, we say that $f_0$ {\it represents} $f$.
If $X={\rm Spec}\, A$, we denote by $\cO_{I,\cP}(X)$ the ring
of families (\ref{fma}) that can be analytically continued along $I$.
\qed
\end{definition}

\medskip

{\it From this point on, the $\delta_{\cP}$-ring $A_0$ under consideration
will be just $\bZ_{(\cP)}$. Cf. the notation in \S \ref{cou}.}
\medskip

Given a $\bZ_{(\cP)}$-point $P: {\rm Spec}\, \bZ_{(\cP)} \ra X$, by
the universality property, we obtain a unique lift to a point
$P^r: {\rm Spec}\, \bZ_{(\cP)} \ra \cJ^r_{\cP}(X)$ that is compatibly with the
action of $\delta_{\cP}$.

\begin{definition}
\label{margica}
Let $X$ be an affine scheme of finite type over $\bZ_{(\cP)}$, $P$ be a
$\bZ_{(\cP)}$-point $P:{\rm Spec}\, \bZ_{(\cP)} \ra X$, and
$P^r:{\rm Spec}\, \bZ_{(\cP)} \ra \cJ^r_{\cP}(X)$ be its unique lift
compatible with the action of $\delta_{\cP}$. By abuse of notation,
we also denote by $P^r \subset \cO(\cJ^r_{\cP}(X))$ the ideal of the image
of $P^r$. By a {\it $\delta_{\cP}$-function
on $X$  that is analytically continued along $P$}  we mean
a family
\begin{equation}
\label{famm}f=(f_k) \in \prod_{k=1}^d
\cO(\cJ^r_{\cP}(X))\whk
\end{equation}
that can be analytically continued along $P^r$.
We denote by $\cO^r_{P,\cP}(X)$ the ring of all $\delta_{\cP}$-functions on
$X$ that are analytically continued  along $P$.
\qed
\end{definition}

Let us recall that by Definition \ref{ideeea}, if $f=(f_k) \in \prod_{k=1}^d
\cO(\cJ^r_{\cP}(X))\whk $ is
analytically continued along $P$ there exists
an element
$$
f_0 \in \cO(\cJ^r_{\cP}(X))^{\widehat{P^r}}
$$
that represents $f$ such that the
images of $f_0$ and $f_k$ in
$$\cO(\cJ^r_{\cP}(X))^{\widehat{(p_k,P^r)}}$$ coincide
for each $k=1,\ldots ,d$. Thus, the the ring of all $\delta_{\cP}$-functions
on $X$ that are analytically continued  along $P$
is
$$
\cO^r_{P,\cP}(X):=\cO_{P^r,\cP}(\cJ^r_{\cP}(X))\, .
$$

In the next section, we shall encounter natural examples of
families (\ref{famm}) that can be analytically continued along a point
$P$ but not along others, as well as examples of such families
that can be analytically continued along {\it any} point.

Under very general hypotheses on $X$ and $P$ that suffice for our
applications, this analytical continuation concept admits a practical
description that we now explain.

\begin{definition}
\label{vocc}
Let $X$ be a smooth affine scheme over $\bZ_{(\cP)}$. A $\bZ_{(\cP)}$-point
$P:{\rm Spec}\, \bZ_{(\cP)}\ra X$ of $X$ is called {\it uniform} if there
exists an \'{e}tale map of $\bZ_{(\cP)}$-algebras $\bZ_{(\cP)}[T] \ra \cO(X)$,
where $T$ is a tuple of indeterminates, such that the ideal of
the image of $P$ in $\cO(X)$ is generated by $T$. We refer to $T$ as {\it uniform coordinates}.
Then $T$ is a regular sequence in $\cO(X)$, and the graded ring associated
to the ideal $(T)$ in $\cO(X)$ is isomorphic to $\bZ_{(\cP)}[T]$;
cf. \cite{mats}, p. 125. We have that
$$
\cO(X)^{\widehat{T}} \simeq \bZ_{(\cP)}[[T]]\, ,
$$
and similarly that
$$
\cO(X)^{\widehat{(p_k,T)}}\simeq \bZ_{p_k}[[T]]\, .
$$
For a general scheme $X$, we say that a $\bZ_{(\cP)}$-point $P$ of $X$ is
{\it uniform} if there exists an affine open set $X_1 \subset X$ that
contains $P$ such that $P$ is uniform in $X_1$.
\qed
\end{definition}

Uniformity is a property that holds ``generically'' in the following sense:
if $X_{\bQ}$ is a smooth scheme over $\bQ$ and $P_{\bQ}$ is a $\bQ$-point of
$X_{\bQ}$, then there exists an integer $N \geq 1$ and models
$X_{\bZ[1/N]}$, $P_{\bZ[1/N]}$ over $\bZ[1/N]$ such that, for any set of
primes $\cP$ that are coprime to $N$, the induced point $P_{\bZ_{(\cP)}}$ of
$X_{\bZ_{(\cP)}}$ is uniform.

\begin{remark}
Let $X$ be a smooth affine scheme with connected geometric fibers over
$\bZ_{(\cP)}$ provided with its tautological cover, and let
$P$ be a uniform point with uniform coordinates $T$.
We consider the corresponding \'{e}tale map $\bZ_{(\cP)}[T]\ra \cO(X)$,
 and let $t$ be the tuple of indeterminates
$(\delta_{\cP}^iT)_{i \leq r}$. By \cite{borger}, Proposition 8.2, the map
$$
\bZ_{(\cP)}[t]\ra \cO(\cJ^r_{\cP}(X))
$$
is \'{e}tale, and in particular, the $r$-jet space $\cJ^r_{\cP}(X)$ is smooth
over $\bZ_{(\cP)}$. Once again, we let $P^r$ be the canonical lift of $P$. Then
$P^r$ is uniform with uniform coordinates $t$. By Definition \ref{vocc}, we have
that
$$
\begin{array}{lcl}
\cO(\cJ^r(X))^{\widehat{t}} & \simeq & \bZ_{(\cP)}[[t]]\, ,\\
\cO(\cJ^r(X))^{\widehat{(p_k,t)}} & \simeq & \bZ_{p_k}[[t]]\, .
\end{array}
$$
Hence, by Definition \ref{margica}, a family
\begin{equation}
\label{fuhh}
f=(f_k) \in
 \prod_{k=1}^d \cO({\mathcal J}^r_{\cP}(X))\whk
\end{equation}
is a $\delta_{\cP}$-function on $X$ that is analytically continued along $P$
if there exists $f_0 \in \bZ_{(\cP)}[[t]]$ such that
the images of $f_k$ and $f_0$ in $\bZ_{p_k}[[t]]$ coincide for
each $k=1,\ldots ,d$.
Notice that such an $f_0$ is uniquely determined by $f$.
\qed
\end{remark}

\begin{definition}
Let $X$ be a smooth quasi-projective scheme over $\bZ_{(\cP)}$ with
geometrically connected fibers, and let $P$
be a uniform point in some affine open set $X_1$.
Let us denote by $\cO_{P,\cP}^r(X)$ the preimage of $\cO_{P,\cP}^r(X_1)$ via
the restriction map
$$
\prod_{k=1}^d \cO({\mathcal J}^r_{\cP}(X)\whk  \ra \prod_{k=1}^d
\cO({\mathcal J}^r_{\cP}(X_1)\whk) =
\prod_{k=1}^d \cO({\mathcal J}^r_{\cP}(X_1))\whk \, .
$$
Elements of this ring $\cO_{P,\cP}^r(X)$ are referred to as
{\it $\delta_{\cP}$-functions of order $r$} on $X$ that are analytically
continued along $P$. We define the
ring of {\it $\delta_{\cP}$-functions on $X$ that are analytically continued
along $P$} by
$$
\cO_{P,\cP}^{\infty}(X):=\lim_{\ra}\cO_{P,\cP}^r(X)\, ,
$$
with its natural $\delta_{\cP}$-ring structure.
\qed
\end{definition}

\begin{example}
Let $X=\bA^n={\rm Spec}\, \bZ_{(\cP)}[T]$ be the affine space of dimension $n$
over $\bZ_{(\cP)}$, where $T$ is an $n$-tuple of indeterminates. We let
$P$ be the zero section, and let $t$ be the tuple
of indeterminates $(\delta_{\cP}^iT)_{i\leq r}$. If
$v_{p_k}$ denotes the $p_k$-adic valuation, then
$$
\cO^r_{P,\cP}(\bA^n)\simeq \{\sum a_jt^j\in \bZ_{(\cP)}[[t]]\, : \;
\text{$\lim_{|j|\rightarrow \infty}v_{p_k}(a_j)\ra \infty$ for each $k$}\}\, .
$$
\qed
\end{example}

\begin{remark}
Let $u:X \ra Y$ be a morphism of quasi-projective schemes over $\bZ_{(\cP)}$.
Assume that $X$ and $Y$ are smooth with connected geometric fibers. Let $P$
be a uniform point of $X$ such that $u(P)$ is a uniform point of $Y$.
Then there is a naturally induced homomorphism
$$
u^*:\cO^r_{u(P),\cP}(Y) \ra \cO^r_{P,\cP}(X)\, .
$$
\qed
\end{remark}

\begin{remark}
\label{maghere}
By Lemma \ref{YO} it follows that if $X$ is a smooth quasi-projective scheme
with connected geometric fibers over $\bZ_{(\cP)}$, and $P$ is a uniform
point of $X$, then the induced homomorphisms
$$
\cO(\cJ^r_{\cP}(X)\whk) \ra \bZ_{p_k}[[\delta_{\cP}^iT :\, i\leq r]]
$$
are injective. In particular, the homomorphism
$$
\cO^r_{P,\cP}(X) \ra \bZ_{(\cP)}[[\delta_{\cP}^iT : \, i\leq r]]
$$
is injective.
\qed
\end{remark}

\begin{definition}
\label{sweetadele}
Let $X$ be a smooth quasi-projective scheme over $\bZ_{(\cP)}$ with
geometrically connected fibers, and let $P$
be a uniform point of $X$. We consider
a $\delta_{\cP}$-function of order $r$ analytically continued along $P$,
$f =(f_{k})\in \cO^r_{P,\cP}(X)$. If $A$ is a $\delta_{\cP}$-ring over
$\bZ_{(\cP)}$, then there is an induced map of sets
$$
f_{A}:X(A_{\cP})\ra A_{\cP}\, ,
$$
where $A_{\cP}=\prod A\whk$ (cf. the notation in \S \ref{cou}), defined as
follows. Let $Q=(Q_{k})\in X(A_{\cP})=\prod_{k=1}^d X(A\whk)$ be a point.
Since each $A\whk$ is a $\delta_{\cP}$-ring
(cf. Remark \ref{comaredel}), by the universality property of
$\delta_{\cP}$-jet spaces the morphism
$Q_{k}:{\rm Spec}\, A\whk \ra X$ lifts to a morphism
$Q_{k}^r: {\rm Spec}\, A\whk \ra \cJ^r_{c,\cP}(X)$, where $c$ is any
principal covering. Since $A\whk$ is $p_k$-adically complete, we obtain an
induced morphism $\widehat{Q}_{k}^r:{\rm Spf}\, A\whk \ra \cJ^r_{c,\cP}(X)
\whk$, and therefore a homomorphism
$\cO(\cJ^r_{\cP}(X)\whk)\ra A\whk$. We denote the image of $f_{k}$ under this
last homomorphism by $f_{k}(Q_{k})$. Then
$$
f_{A}(Q):=(f_{k}(Q_{k})) \in A_{\cP}\, .
$$
\qed
\end{definition}

Using the map $f_{A}$ of Definition \ref{sweetadele} above,
we may speak of the {\it space of solutions}, $f_{A}^{-1}(0)$.

Notice that the  definition of $f_A$ is functorial in $A$. For if
$A \ra B$ is a morphism of $\delta_{\cP}$-rings, then the corresponding
maps $f_A$ and $f_B$ are compatible with each other.

\subsection{$\delta_{\cP}$-characters}
Let $G$ be a smooth group scheme over $\bZ_{(\cP)}$ with multiplication
$\mu:G \times G \ra G$. Here, $\times=\times_{\bZ_{(\cP)}}$.
If $G$ is affine, by the universality property we have that
${\mathcal J}_{\cP}^r(G)$ is a group scheme over $\bZ_{(\cP)}$.
In the non-affine case this is not a priori true because our definitions
are covering dependent. However, we attach to this general situation a formal
group law as follows.

Let $Z$ be  the identity ${\rm Spec}\, \bZ_{(\cP)} \ra G$, and assume $Z$ is
uniform, with uniform coordinates $T$. Then $Z \times Z \subset G \times G$ is
a uniform point with uniform coordinates $T_1,T_2$ induced by $T$.
We have an induced homomorphism $\bZ_{(\cP)}[[T]] \ra \bZ_{(\cP)}[[T_1,T_2]]$
that sends the ideal $(T)$ into the ideal $(T_1,T_2)$. B
y the universality property applied to the restriction
$\bZ_{(\cP)}[T]\ra \bZ_{(\cP)}[[T_1,T_2]]$, we obtain morphisms
$$
\bZ_{(\cP)}[\delta_{\cP}^i T : \, i\leq r]\ra \bZ_{(\cP)}[[\delta_{\cP}^iT_1,
\delta_{\cP}^iT_2 : \, i\leq r]]
$$
that send the ideal generated by the variables into the corresponding ideal
generated by the variables.
Thus, we have an induced morphism
$$
\bZ_{(\cP)}[[\delta_{\cP}^i T: \, i\leq r]]\ra \bZ_{(\cP)}
[[\delta_{\cP}^iT_1,\delta_{\cP}^iT_2: \, i\leq r]]\, .
$$
We call ${\mathcal G}^r$ the image of the variables
$\{\delta_{\cP}^i T: \, i\leq r\}$ under this
last homomorphism.
Then the tuple ${\mathcal G}^r$ is a formal group law over $\bZ_{(\cP)}$.

\begin{definition}
Let $G$ be a quasi-projective smooth group scheme over $\bZ_{(\cP)}$ with
geometrically connected fibers. Let us assume that the identity is a uniform
point. Then there are homomorphisms
$$
\mu^*,pr_1^*,pr_2^*:\cO^r_{Z,\cP}(G) \ra \cO^r_{Z \times Z,\cP}(G \times G)
$$
induced by the product $\mu$ and the two projections. We say that
a $\delta_{\cP}$-function $f \in \cO^r_{Z,\cP}(G)$ of order $r$ on $G$ that
is analytically continued along $Z$ is a
{\it $\delta_{\cP}$-character of order $r$ on $G$} if
$$
\mu^* f=pr_1^* f+pr_2^* f\, .
$$
We denote by $\bX^r_{\cP}(G)$ the group of $\delta_{\cP}$-characters of
order $r$ on $G$. We define
the {\it group of $\delta_{\cP}$-characters} on $G$ to be
$$
\bX^{\infty}_{\cP}(G):=\lim_{\ra} \bX^r_{\cP}(G)\, .
$$
\qed
\end{definition}

By Remark \ref{maghere}, the condition that $f$ be a
$\delta_{\cP}$-character of order $r$ on $G$ is that there exists
$$
f_0 \in \bZ_{(\cP)}[[\delta_{\cP}^iT: \, i\leq r]]
$$
that represents $f$ such that
\begin{equation}
\label{cocca}
f_0({\mathcal G}^r(T_1,T_2))=f_0(T_1)+f_0(T_2)\, .
\end{equation}
Here, $\cG^r$ is the corresponding formal group law, and
$f_0(T)$ stands for $f_0(\ldots ,\delta_{\cP}^i T,\ldots )$.

The group $\bX^r_{\cP}(G)$ of $\delta_{\cP}$-characters of order $r$ on $G$ is
a subgroup of the additive group of the ring
$\cO^r_{Z,\cP}(G)$. The group $\bX^{\infty}_{\cP}(G)$ of
$\delta_{\cP}$-characters on $G$
is a subgroup of $\cO_{Z,G}^{\infty}(G)$.

For any $\delta_{\cP}$-character $\psi \in \bX_{\cP}^r(G)$ and any
$\delta_{\cP}$-ring $A$, there is an induced group homomorphism
\begin{equation}
\label{indho}
\psi_{A}:G(A_{\cP})\ra A_{\cP}\, ,
\end{equation}
where $A_{\cP}$ is viewed as a group with respect to addition.
We may therefore speak of the {\it group of solutions}
${\rm Ker}\, \psi_{A}$.

Once again, the mapping $A \mapsto \psi_A$ is functorial in $A$.

\begin{definition}
Let $P$ be a $\bZ_{(\cP)}$-point of $G$.
We say that a $\delta_{\cP}$-character $\psi \in \bX_{\cP}^r(G)$ can
be {\it analytically continued along $P$}
if $\psi \in \cO^r_{P,\cP}(G)$.
\qed
\end{definition}

By the mere definition, any $\delta_{\cP}$-character can be analytically
continued along the identity $Z$. Later on we will address the question of
when a $\delta_{\cP}$-character can be analytically continued along points
other than $Z$.

\section{Main results}
\setcounter{theorem}{0}

In this section we determine all the  $\delta_{\cP}$-characters on the
additive group $\bG_a$, the multiplicative group $\bG_m$,
and elliptic curves $E$ over $\bZ_{(\cP)}$.

\subsection{The additive group}
We consider the additive group scheme over $\bZ_{(\cP)}$,
$$
\bG_a:={\rm Spec}\, \bZ_{(\cP)}[x]\, .
$$
The zero section is uniform, with uniform coordinate
$T=x$. We have
$$
\cO({\mathcal J}_{\cP}^r(\bG_a)) = \bZ_{(\cP)}
[\delta_{\cP}^i x : \, i\leq r]\, .
$$
Let us consider the polynomial ring
$$
\bZ_{(\cP)}[\phi_{\cP}]:=
\bZ_{(\cP)}[\phi_{p_l}: \, l \in I]=
\sum_{n \in {\mathcal N}} \bZ_{(\cP)} \phi_n \, ,
$$
where $I:=\{1,\ldots ,d\}$, ${\mathcal N}$ is the monoid of the natural
numbers generated by $\cP$, and the $\phi_{p_l}$s are commuting variables.
For $i=(i_1,\ldots ,i_d)\in \bZ_+^d$ we set
$\cP^i=p_1^{i_1} \ldots p_d^{i_d}$. If $n=\cP^i$ we set
$\phi_n=\phi_{\cP}^i=\phi_{p_1}^{i_1}\ldots \phi_{p_d}^{i_d}$.
We will view the ring $\bZ_{(\cP)}[\phi_{\cP}]$ as the
{\it ring of symbols}; cf. \cite{pde}.
Let $r \in \bZ_+^d$ and
$$
\bar{\psi}:=\sum_{n |\cP^r} c_n \phi_n\in \bZ_{(\cP)}[\phi_{\cP}]\, .
$$
We may consider the element
$\psi=\bar{\psi}x \in \cO({\mathcal J}^r(\bG_a))$ and identify it
with an element
$$
\psi \in \prod_{k=1}^d \cO({\mathcal J}_{\cP}^r(\bG_a))\whk \, .
$$
via the diagonal embedding. Then $\bar{\psi}$
 clearly defines a $\delta_{\cP}$-character on $\bG_a$. We will usually
identify $\psi$ and $\bar{\psi}$.

As we see now, the following result is easy to prove for the
additive group. Later on, we shall prove a less elementary analogue for
the multiplicative group $\bG_m$, and for elliptic curves.

\begin{theorem}
\label{soap}
Let $\psi$ be a  $\delta_{\cP}$-character on $\bG_a$ of order $r$.
Then
\begin{enumerate}
\item $\psi$ can be uniquely written as $\psi=\bar{\psi} x$ with
$$\bar{\psi}=(\sum_{n|\cP^r} c_n \phi_n)x,\ \ c_n \in \bZ_{(\cP)}.$$
\item $\psi$ can be analytically continued along  any $\bZ_{(\cP)}$-point
$P$ of $\bG_a$.
\end{enumerate}
\end{theorem}

\begin{corollary}
The group of $\delta_{\cP}$-characters $\bX^{\infty}_{\cP}(\bG_a)$
is a free $\bZ_{(\cP)}[\phi_{\cP}]$-module of rank one with basis $x$.
\end{corollary}

{\it Proof of Theorem \ref{soap}}.
Let $\psi \in \bX^{r}_{\cP}(\bG_a)$ be represented by a series
$\psi_0 \in \bZ_{(\cP)}[[\delta_{\cP}^i x: \, i\leq r]]$. We view $\psi_0$
as an element of
$$
\bQ[[\delta_{\cP}^i x: \, i\leq r]]=\bQ[[\phi_{\cP}^i x: \, i\leq r]]\, .
$$
We have that
$$
\begin{array}{rcl}
\psi_0(\ldots ,\phi_{\cP}^i x_1+\phi_{\cP}^i x_2,\ldots ) &= &
\psi_0(\ldots ,\phi_{\cP}^i (x_1+ x_2), \ldots ) \\
 & = & \psi_0(\ldots ,\phi_{\cP}^i x_1,\ldots )+\psi_0(\ldots ,
\phi_{\cP}^i x_2,\ldots ) \, .
\end{array}
$$
Then it follows that $\psi_0=\sum c_n \phi_n x$ with $c_n \in \bQ$.

Now notice that
$$
\phi_n x\equiv x^n\; {\rm mod}\, (\delta_{\cP}^i x: \, i\leq r)
$$
in the ring $\bQ[\delta_{\cP}^ix : \, i\leq r]$. We get that
$$
(\psi_0)_{|\delta_{\cP}^i x=0; \, 0 \neq i\leq r}=\sum c_n x^n\, .
$$
It follows that $c_n \in \bZ_{(\cP)}$, and the first assertion is proved.

The second assertion is obvious.
\qed

\subsection{The multiplicative group}
We now consider the multiplicative group scheme over $\bZ_{(\cP)}$,
$$
\bG_m:={\rm Spec}\, \bZ_{(\cP)}[x,x^{-1}]\, .
$$
The zero section is uniform, with uniform coordinate $T=x-1$.

We recall that the formal group law ${\mathcal G}^0$ corresponding to $T$
is
$$
{\mathcal G}^0(T_1,T_2)=T_1+T_2+T_1T_2\, ,
$$
and the logarithm of this formal group law is given by the series
$$
l_{\bG_m}(T)=\sum_{n=1}^{\infty} (-1)^{n-1}\frac{T^n}{n}\in \bQ[[T]]\, .
$$
By Lemma \ref{dorinpiept}, we have that
$$
\cO({\mathcal J}_{\cP}^r(\bG_m))  =  \bZ_{(\cP)}\left[\delta_{\cP}^i x,
\frac{1}{\phi_{\cP}^i(x)}: \, i\leq r\right]\, .
$$
For each $k$, we consider the series
$$
\psi_{p_k}=\psi_{p_k}^1:=\sum_{n=1}^{\infty}(-1)^{n-1}
\frac{p_k^{n-1}}{n}\left(\frac{\delta_{p_k} x}{x^{p_k}} \right)^n
\in \cO({\mathcal J}_{\cP}^{e_k}(\bG_m))\whk\, ,
$$
and the element
$$
\bar{\psi}_{\bar{p}_k} :=  \prod_{l \in I_k} \left( 1-
\frac{\phi_{p_l}}{p_l}\right)=\sum_{n \in {\mathcal N}_k}
\frac{\mu(n)}{n}\phi_n \in \bZ_{(p_k)}[\phi_{p_l}: \, l \in I_k]\, .
$$
Here $I_k:=\{1,\ldots ,d\}\backslash \{k\}$,
${\mathcal N}_k$ is the monoid of the natural numbers generated by
$\bar{p}_k:=\cP\backslash \{p_k\}$, $\mu$ is the Mobius function
and the $\phi_{p_l}$s are commuting variables.
We then may consider the family
$$
(\bar{\psi}_{\bar{p}_k} \psi_{p_k}) \in \prod_{k=1}^d
\cO({\mathcal J}_{\cP}^e(\bG_m))\whk\, ,
$$
where $e=e_1+\cdots +e_d=(1,\ldots ,1)$.

\begin{theorem}
\label{Gmtheorem}
The family $(\bar{\psi}_{\bar{p}_k} \psi_{p_k})$ is a
$\delta_{\cP}$-character of order $e$ on $\bG_m$.
\end{theorem}

We denote this $\delta_{\cP}$-character by
\begin{equation}
\label{thecharm}
\psi^e_m \in \bX^{e}_{\cP}(\bG_m)\, .
\end{equation}

{\it Proof.} We notice that
\begin{equation}
\label{bufff}
\delta_{p_k}(1+T) = \delta_{p_k}T+C_{p_k}(1,T) \in (T,\delta_{p_k}T)
\subset \bZ_{p_k}[[T,\delta_{p_k}T]]\, .
\end{equation}
Then the image of $\bar{\psi}_{\bar{p}_k} \psi_{p_k}$ in
$\bQ_{p_k}[[\delta_{\cP}^i T: \, i\leq e]]$ is equal
the following series (that, by (\ref{bufff}), is convergent in the
topology given by the maximal ideal of this ring):
$$
\begin{array}{ll}
\bar{\psi}_{\bar{p}_k} \left( \sum_{n=1}^{\infty}  \frac{(-p_k)^{n-1}}{n}
\left(
\frac{\delta_{p_k}(1+T)}{(1+T)^{p_k}}\right)^n\right)
& = {\displaystyle
 \frac{1}{p_k} \bar{\psi}_{\bar{p}_k}  \left( \sum_{n=1}^{\infty}
\frac{(-1)^{n-1}}{n} \left(
\frac{\phi_{p_k}(1+T)}{(1+T)^{p_k}}-1\right)^n\right)} \vspace{1mm} \\
 & ={\displaystyle \frac{1}{p_k}\bar{\psi}_{\bar{p}_k}  l_{\bG_m}
\left(\frac{\phi_{p_k}(1+T)}{(1+T)^{p_k}}-1\right)} \vspace{1mm} \\
 & ={\displaystyle
\frac{1}{p_k}\bar{\psi}_{\bar{p}_k} (\phi_{p_k}-p_k)l_{\bG_m}(T)}\vspace{1mm}\\
& ={\displaystyle
-\left(\prod_{l=1}^d(1-\frac{\phi_{p_l}}{p_l})\right) l_{\bG_m}(T)}\, .
\end{array}
$$
We call $\psi^e_0$ this series.

We have then that $\psi^e_0$ has coefficients in $\bQ\cap \bZ_{p_k}=
\bZ_{(p_k)}$, and is the same for all $k=1,\ldots ,d$. Hence,
$\psi^e_0$ has coefficients in $\bZ_{(\cP)}$ and represents the family
$(\bar{\psi}_{\bar{p}_k} \psi_{p_k})$.
Further, $\psi^e_0$ satisfies the condition (\ref{cocca}) because
$$
\begin{array}{lll}
\psi^e_0({\mathcal G}^e(T_1,T_2)) &
 = &  -\left(\prod_{l=1}^d(1-\frac{\phi_{p_l}}{p_l})\right)
\left(l_{\bG_m}(T)({\mathcal G}^0(T_1,T_2)\right) \\
 &  = &  -\left(\prod_{l=1}^d(1-\frac{\phi_{p_l}}{p_l})\right)
\left(l_{\bG_m}(T_1)+l_{\bG_m}(T_2)\right) \vspace{1mm} \\
 & = &  {\displaystyle \psi^e_0 (T_1)+
 \psi^e_0 (T_2)} \, .
\end{array}
$$
\qed

We now prove that  $\psi^e_m$ (cf. (\ref{thecharm}))
generates, in an appropriate sense, the space of all $\delta_{\cP}$-characters
of $\bG_m$. We also determine which $\delta_{\cP}$-characters $\psi$ can
be analytically continued along any given  $\bZ_{(\cP)}$-point $P$ of $\bG_m$.

\begin{theorem}
\label{ungm}
Let $\psi$ be a $\delta_{\cP}$-character of order $r$ on $\bG_m$. Then
\begin{enumerate}
\item $\psi$ can be uniquely written as
$$
\psi=(\sum_{n} \rho_n \phi_n) \psi_m^e\, ,\quad \rho_n \in \bZ_{(\cP)}\, .
$$
\item $\psi$ can be analytically continued along a $\bZ_{(\cP)}$-point
$P$ of $\bG_m$ if, and only if, either $P$ is torsion
or $\sum_n \rho_n=0$.
\end{enumerate}
\end{theorem}

Consider the {\it augmentation ideal}
$$
\bZ_{(\cP)}[\phi_{\cP}]^+:=\{\sum_n \rho_n \phi_n \in \bZ_{(\cP)}[\phi_{\cP}]:
\, \sum_n \rho_n=0\}\, .
$$

\begin{corollary}
The group of $\delta_{\cP}$-characters $\bX^{\infty}_{\cP}(\bG_m)$
is a free $\bZ_{(\cP)}[\phi_{\cP}]$-module of rank one with basis $\psi_m^e$.
The group of $\delta_{\cP}$-characters in $\bX^{\infty}_{\cP}(\bG_m)$
that can be analytically continued along a given non-torsion point $P$ of
$\bG_m$ is isomorphic with the augmentation ideal $\bZ_{(\cP)}[\phi_{\cP}]^+$
as a $\bZ_{(\cP)}[\phi_{\cP}]$-module. This group is the same for all
non-torsion $P$s.
\end{corollary}

{\it Proof of Theorem \ref{ungm}}.
Let $\psi_0$ be the series representing $\psi$.
We first prove the following:

\begin{enumerate}
\item[Claim 1.] There is an equality
$$\psi_0=(\sum_{n|\cP^r} c_n \phi_n)l(T)$$
in $\bQ[[\delta_{\cP}^iT: \, i\leq r]]$ with $c_n \in \bQ$.
\medskip

Indeed let $e(T) \in \bQ[[T]]$ be the compositional inverse of $l(T)$.
Then the series
$$
\Theta(\ldots ,\delta_{\cP}^i(T),\ldots ):=
\psi_0(\ldots ,\delta_{\cP}^i(e(T)),\ldots )
$$
satisfies
$$
\Theta(\ldots ,\delta_{\cP}^i(T_1+T_2),\ldots )=\Theta(\ldots ,
\delta_{\cP}^i T_1,\ldots )+
\Theta(\ldots ,\delta_{\cP}^i T_2,\ldots )\, .
$$
As in the proof of Proposition \ref{soap},
we conclude that $\Theta=\sum c_n \phi_n(T)$, with $c_n \in \bQ$,
which finishes the proof of the claim.
\end{enumerate}

By Claim 1, by setting $\delta_{\cP}^i T=0$ for $i\neq 0$, we obtain
$$
(\sum c_n \phi_n) \star l(T) \in \bZ_{(\cP)}[[T]]\, ,
$$
where $\phi_n \star T:=T^n$.
\medskip

\begin{enumerate}
\item[Claim 2.] If a polynomial $\Lambda=\sum \lambda_n \phi_n \in
\bQ[\phi_{p_1},\ldots ,\phi_{p_s}]$
satisfies
$$
\Lambda \star l(T) \in \bZ_{p_k}[[T]]\otimes \bQ
$$
for some $k\in \{1,\ldots ,s\}$, then $\Lambda$ is divisible in
the ring  $\bQ[\phi_{p_1},\ldots ,\phi_{p_s}]$ by $\phi_{p_k}-p_k$.
\medskip

Indeed, let us divide $\Lambda$ by
$\phi_{p_k}-p_k$ in $\bQ[\phi_{p_1},\ldots ,\phi_{p_s}]$ to obtain
$$
\Lambda=(\sum \alpha_n \phi_n)(\phi_{p_k}-p_k)+\sum \beta_n \phi_n\, ,
$$
where $\alpha_n,\beta_n \in \bQ$ and $\beta_n=0$ if $p_k|n$. We analyze
the remainder term, and prove that it is identically zero by showing that
$\beta_n=0$ for all $n$.

Since $(\phi_{p_k}-p_k) \star l(T) \in \bZ_{p_k}[[T]]$, it follows that
$(\sum \beta_n \phi_n)\star l(T) \in \bZ_{p_k}[[T]]\otimes \bQ$.
We may assume $(\sum \beta_n \phi_n)\star l(T) \in \bZ_{p_k}[[T]]$.
We have that
\begin{equation}
\label{iese}
(\sum \beta_n \phi_n)\star l(T) =\sum_n \sum_m (-1)^{m-1} \beta_n
\frac{T^{nm}}{m}\, .
\end{equation}

Let us fix integers $n',\nu \geq 1$. By looking at the coefficient of
$T^{n'p_k^{\nu}}$ in (\ref{iese}), we obtain that
$$
-\sum_{n|n'} (-1)^{n'/n}\frac{n \beta_n}{n'p_k^{\nu}} \in \bZ_{p_k}
$$
because the equality $nm=n'p_k^{\nu}$ with $n \not\equiv 0$ mod $p_k$
implies $m=\mu p_k^{\nu}$,
 $\mu \in \bZ$, $n\mu=n'$; so $n|n'$ and $m=\frac{n'}{n} p_k^{\nu}$.
Since $n$ is odd for $\beta_n \neq 0$, it follows that
$$
\sum_{n|n'} n \beta_n\in p_k^{\nu}\bZ_{p_k}\, ,
$$
and since this is true for all $\nu$, we obtain that
$$
\sum_{n|n'} n \beta_n=0\, .
$$
But this is true for all $n'$. By the  Mobius' inversion formula, we
conclude that $\beta_n=0$ for all $n$. This completes the proof of
Claim 2.
\end{enumerate}

By Claim 2, it follows that
$$
\sum c_n\phi_n=(\sum \rho_n \phi_n)\prod_{k=1}^d \left( 1-
\frac{\phi_{p_k}}{p_k}\right)
$$
for some $\rho_n \in \bQ$. By induction on $n$, it is then easy to check
that $\rho_n \in \bZ_{(\cP)}$ for all $n$. This completes the proof of the
first assertion (1) in the statement.

In order to prove assertion (2), let $\tau:\bG_m \ra \bG_m$ be the
translation defined by the inverse of $P$, and let $\tau^*$ be the
automorphism defined by $\tau$ on the various rings of functions.
Since  $\psi\in \cO_{Z,\cP}^r(\bG_m)$, we have that
 $\tau^* \psi \in \cO^r_{P,\cP}(\bG_m)$. But if $\psi=(\psi_k)$
then $\tau^* \psi_k= \psi_k-\psi_k(P_k)$ for all $k$; cf. the notation in
Definition \ref{sweetadele}. Now, if $P$ is torsion
or if $\sum_n \rho_n=0$, it is then clear that
$\psi_k(P_k)=0$.
So $\tau^* \psi_k=\psi_k$, hence $\psi \in \cO^r_{P,\cP}(\bG_m)$.

Conversely, let $P$ be non-torsion and given by a number
$a \in \bZ_{(\cP)}^{\times}$. Let $\sum_n \rho_n\neq 0$, and assume
$\psi\in \cO^r_{P,\cP}(\bG_m)$. We derive a contradiction.
For let $p=p_1$ and $b:=a^{p-1} \in 1+p\bZ_{(p)}$, so $b \neq 1$.
By Mahler's $p$-adic analogue of the Hermite-Lindemann theorem
\cite{mahler,ber}, we have that $\log{(b)} \not\in \bQ$
(where
$\log:1+p\bZ_p\ra p\bZ_p$ is the $p$-adic logarithm).
Since $\tau^* \psi \in \cO^r_{P,\cP}(\bG_m)$ and $\psi \in
\cO^r_{P,\cP}(\bG_m)$, it follows that $\tau^*\psi-\psi
\in \cO^r_{P,\cP}(\bG_m)$. But $\tau^*\psi-\psi=(-\psi_k(a))$ so, in
particular,  $\psi_1(b) \in \bQ$.
But
$$\psi_1(b)=-(\sum_n \rho_n) \cdot \prod_{l=1}^d
\left(1-\frac{1}{p_l}\right) \cdot \log{ b}\, ;
$$
cf. the proof of Theorem \ref{Gmtheorem}. Since $\sum_n\rho_n\neq 0$,
it follows that $\log{b}\in \bQ$, reaching the desired contradiction.
\qed

The next result computes the group of solutions of the
$\delta_{\cP}$-character $\psi^e_m$ (\ref{thecharm}).

\begin{theorem}
\label{oastept}
Let $A$ be the $\delta_{\cP}$-ring $\bZ_{(\cP)}[\zeta_m]$ in
Example \ref{basicex}, and let $\psi^e_{m,A}:\bG_m(A_{\cP})=
A_{\cP}^{\times} \ra A_{\cP}$ be the homomorphism (\ref{indho}) induced
by $\psi^e_m$. Then
$$
{\rm Ker}\, \psi^e_{m,A}=(A_{\cP}^{\times})_{tors}\, .
$$
\end{theorem}

In the statement above, we use the notation $\Gamma_{tors}$ to denote the
torsion group of an abelian group $\Gamma$.

{\it Proof.} The non-trivial inclusion is ``$\subset$.'' Let us
take $Q=(Q_{k}) \in {\rm Ker}\, \psi^e_{m,A}$
so
\begin{equation}
\label{pott}
\bar{\psi}_{\bar{p}_k}(\psi_{p_k}(Q_{k}))=0
\end{equation}
for all $k$. Here $Q_{k}\in A\whk=A^{\widehat{P_{k1}}}
\times A^{\widehat{P_{k2}}} \times \cdots $, where
$p_k=P_{k1}P_{k2}\cdots $ is the prime decomposition of $p_kA$. In order
to show that $Q$ is torsion, we may replace $Q$ by any of its
powers. So we may assume that $Q_{k}\in 1+p_k A\whk$ for all $k$.
Then (\ref{pott}) gives
$$
(\prod_{l=1}^d (\phi_{p_l}-p_l))l_{\bG_m}(Q_{k}-1)=0\, .
$$
We claim that the map
$$
\begin{array}{rcl}
A\whk & \rightarrow & A\whk \\
\beta & \mapsto & (\phi_{p_l}-p_l)\beta
\end{array}
$$
is injective for all $k$, $l$. Using this claim, we conclude that
$l_{\bG_m}(Q_{k}-1)=0$, which implies that $Q_{k}=1$ by the injectivity
of $l_{\bG_m}:p_kA\whk \ra p_kA\whk$, and finishes the proof of the
Theorem.

In order to prove the claim, let us assume that $(\phi_{p_l}-p_l)\beta=0$. We
also have that $\phi_{p_l}^M \beta=\beta$, $M:=[\bQ(\zeta_m):\bQ]$. Since
the polynomials
$\phi_{p_l}-p_l, \phi_{p_l}^M-1 \in \bQ[\phi_{p_l}]$ are coprime,
it follows that $\beta=0$, as desired.
\qed

\subsection{Elliptic curves}

Let $E$ be an elliptic curve over $\bQ$. We assume that $E$ has minimal model
over $\bZ$ given by
$$
y^2+c_1xy+c_3y=x^3+c_2x^2+c_4x+c_6\, .
$$
We assume further that the discriminant of $E$ is not divisible by any of
the primes in $\cP$, and view $E$ as an elliptic curve (smooth projective
group scheme) over $\bZ_{(\cP)}$. Then the zero section is uniform, with
uniform coordinate
 $T=\frac{x}{2y}$. Let ${\mathcal G}^0\in \bZ_{(\cP)}[[T_1,T_2]]$ be the
formal group law attached to $E$ with respect to $T$, and let
$l_E \in \bQ[[T]]$ be the logarithm of ${\mathcal G}^0$.
So
$$
l_E(T)=\sum b_n \frac{T^n}{n}\, ,
$$
where
$$
\frac{dx}{2y+c_1 x+c_3}=\sum b_nT^{n-1}dT\, .
$$

Let $a_{p_k} \in \bZ$ be the trace of the $p_k$-power Frobenius on the
reduction mod $p_k$
of $E$. By (\ref{pisacraz}) and \cite{frob}, Theorem 1.10, there exists
$$
\psi_{p_k}=\psi_{p_k}^2 \in \cO({\mathcal J}^2_{\{p_k\}}(E)\whk)
$$
whose image via
$$
\cO({\mathcal J}^2_{\{p_k\}}(E)\whk)\ra \bQ_{p_k}[[T,\delta_{p_k}
T,\delta^2_{p_k}T]]
$$
equals
$$
\frac{1}{p_k}(\phi_{p_k}^2-a_{p_k}\phi_{p_k}+p_k)l_E(T)\, .
$$
On the other hand, we may consider the element
$$
\bar{\psi}_{\bar{p}_k}: = \prod_{l \in I_k} \left(
 1-a_{p_l}\frac{\phi_{p_l}}{p_l}+p_l\left( \frac{\phi_{p_l}}{p_l}\right)^2
 \right)\in \bZ_{
p_k}[\phi_{p_l}: \, l \in I_k]\, ,
$$
where $I_k=\{1,\ldots ,d\}\backslash \{k\}$. Let us still denote
by $\psi_{p_k}$ the image of $\psi_{p_k}$ via the homomorphism
$$
\cO({\mathcal J}_{\{p_k\}}^2(E)\whk)
\ra \cO({\mathcal J}_{\cP}^{2e_k}(E)\whk)
$$
induced by (\ref{pisacraz}).
Then, we may consider the family
$$
(\bar{\psi}_{\bar{p}_k} \psi_{p_k}) \in \prod_{k=1}^d
\cO({\mathcal J}_{\cP}^{2e}(E)\whk)\, .
$$

\begin{theorem}\label{exE}
The family $(\bar{\psi}_{\bar{p}_k} \psi_{p_k})$ is a $\delta_{\cP}$-character
of order $2e$ on $E$.
\end{theorem}

We denote this $\delta_{\cP}$-character by
\begin{equation}
\label{thechare}
\psi^{2e}_E\in \bX^{2e}_{\cP}(E)\, .
\end{equation}

{\it Proof}. As in the proof of
Theorem \ref{Gmtheorem},
the image of $\bar{\psi}_{\bar{p}_k} \psi_{p_k}$ in the ring
$$
\bQ_{p_k}[[\delta_{\cP}^iT\ ;\ i\leq 2e]]
$$
is equal to the series
$$
\psi^{2e}_0:=\left[ \prod_{l=1}^d\left(1-a_{p_l}
\frac{\phi_{p_l}}{p_l}+p_l\left(
\frac{\phi_{p_l}}{p_l}\right)^2\right)\right]l_E(T)\, ,
$$
which is independent of $k$ and has coefficients in $\bZ_{(p_k)}$. Also,
as in loc.cit. we have that
$$
\psi^{2e}_0({\mathcal G}^{2e}(T_1,T_2))=\psi^{2e}_0(T_1)+
\psi^{2e}_0(T_2)\, .
$$
\qed

As in the case of the multiplicative group $\bG_m$, we now have
the following.

\begin{theorem}
\label{unE}
Let us assume that the elliptic curve $E$ has ordinary (good) reduction at
all the primes in $\cP$. Let $\psi$ be a $\delta_{\cP}$-character of order
$r$ on $E$. Then
\begin{enumerate}
\item $\psi$ can be uniquely written as
$$
\psi=(\sum_{n} \rho_n \phi_n) \psi_E^{2e}\, , \quad \rho_n \in \bZ_{(\cP)}\, .
$$
\item $\psi$ can be analytically continued along a $\bZ_{(\cP)}$-point $P$
of $E$ if, and only if, either $P$ is torsion
or $\sum_n \rho_n=0$.
\end{enumerate}
\end{theorem}

\begin{corollary}
 The group of $\delta_{\cP}$-characters $\bX^{\infty}_{\cP}(E)$
is a free $\bZ_{(\cP)}[\phi_{\cP}]$-module of rank one
 with basis $\psi^{2e}_E$.
The group of $\delta_{\cP}$-characters in $\bX^{\infty}_{\cP}(E)$
that can be analytically continued along a given non-torsion point $P$ of $E$
is isomorphic
with the augmentation ideal $\bZ_{(\cP)}[\phi_{\cP}]^+$ as a
 $\bZ_{(\cP)}[\phi_{\cP}]$-module. This
group is the same for all non-torsion $P$s.
\end{corollary}

{\it Proof of Theorem \ref{unE}}. Proceeding as in the proof of
Theorem \ref{ungm}, if $\psi_0$ is the series representing $\psi$, then
the following claim holds:

\begin{enumerate}
\item[Claim 1.] There is an equality
$$\psi_0=(\sum_{n|\cP^r} c_n \phi_n)l_E(T)$$
in $\bQ[[\delta_{\cP}^iT\ ;\ i\leq r]]$ with $c_n \in \bQ$.
\end{enumerate}

By Claim 1,  we obtain that
$$
(\sum c_n \phi_n) \star l_E(T) \in \bZ_{(\cP)}[[T]]\, .
$$

We also have the following:

\begin{enumerate}
\item[Claim 2.]
If a polynomial $\Lambda=\sum \lambda_n \phi_n \in
\bQ[\phi_{p_1},\ldots ,\phi_{p_d}]$
satisfies
$$
\Lambda \star l_E(T) \in \bZ_{p_k}[[T]]\otimes \bQ
$$
for some $k\in \{1,\ldots ,d\}$, then $\Lambda$ is divisible in
the ring  $\bQ[\phi_{p_1},\ldots ,\phi_{p_d}]$ by
$\phi^2_{p_k}-a_{p_k}\phi_{p_k}+p_k$.

Indeed, by the fact that $E$ has ordinary reduction at $p_k$,
 we may consider the unique root
$u_kp_k$ of $x^2-a_{p_k}x+p_k$ in $p_k\bZ_{p_k}$. Since this root is
not in $\bQ$, it is enough
to prove that $\Lambda$ is divisible in $\bQ_{p_k}[\phi_{p_1},\ldots ,
\phi_{p_d}]$ by $\phi_{p_k}-u_kp_k$. Let
$$
L(E,s)=\sum a_n n^{-s}
$$
be the $L$-function of $E$. Let
$$
f_E(T):=\sum a_n \frac{T^n}{n}\, .
$$
As a consequence of the Euler product structure of $L(E,s)$,
it follows that
\begin{equation}
\label{rupt}
(\phi_{p_k}-u_kp_k) \star f_E(T) \in \bZ_{p_k}[[T]]\, ,
\end{equation}
cf. \cite{hazebook}, p. 441. By the Honda-Hill Theorem \cite{haze}, p. 76,
there exists a series $h(T)=T+\cdots \in \bZ[[T]]$ such that
$f_E(T)=l_E(h(T))$. By the Functional Equation Lemma \cite{haze}, p.74,
$$
(\phi_{p_k}-u_kp_k) \star l_E(T) \in \bZ_{p_k}[[T]]\, .
$$
Now, as in the proof of Theorem \ref{ungm}, we divide $\Lambda$ by
$\phi_{p_k}-u_kp_k$ in the ring
$\bQ_{p_k}[\phi_{p_1},\ldots ,\phi_{p_d}]$, to obtain
$$
\Lambda=(\sum \alpha_n \phi_n)(\phi_{p_k}-u_kp_k)+\sum \beta_n \phi_n
$$
with $\alpha_n, \beta_n \in \bQ_{p_k}$, $\beta_n=0$ for $p_k|n$, and proceed
to analyze the remainder.

By (\ref{rupt}), we get
$(\sum \beta_n \phi_n) \star l_E(T) \in \bZ_{p_k}[[T]]\otimes \bQ$.
 We claim that $\beta_n=0$ for all $n$.
We may assume $(\sum \beta_n \phi_n)\star l_E(T) \in \bZ_{p_k}[[T]]$.
We have
\begin{equation}
\label{iesei}
(\sum \beta_n \phi_n)\star l_E(T) =\sum_n \sum_m a_m \beta_n
\frac{T^{nm}}{m}\, .
\end{equation}
Fix integers $n',\nu \geq 1$. As in the proof of Theorem \ref{ungm},
we get
$$
\sum_{n|n'} a_{p_k^{\nu}n'/n} n \beta_n \in p_k^{\nu} \bZ_{p_k}\, .
$$
Assume $n|n'$ and $n' \not\equiv 0$ mod $p_k$. Then
$a_{p_k^{\nu}n'/n}=a_{p_k^{\nu}}a_{n'/n}$.
Since $E$ has ordinary reduction at $p_k$, we have that
$a_{p_k^{\nu}}\not\equiv 0$ mod $p_k$.
We conclude that
$$
\sum_{n|n'}a_{n'/n}n\beta_n \in p_k^{\nu}\bZ_{p_k}\, .
$$
Since $\nu$ is arbitrary, we get
$$
\sum_{n|n'}a_{n'/n}n\beta_n =0\, .
$$
Since this is true for any $n' \not\equiv 0$ mod $p_k$, we get that the
following
product of formal Dirichlet series with coefficients in $\bQ_{p_k}$
vanishes:
$$
(\sum_{n \not\equiv 0 (p_k)} n \beta_n  n^{-s})(\sum_{n\not\equiv
 0 (p_k)}a_n n^{-s})=0\, .
$$
We conclude that $\beta_n=0$ for all $n$, which complets the proof of Claim 2.
\end{enumerate}

>From this point on, we may proceed as in the proof of Theorem \ref{ungm}
to derive assertion (1). Assertion (2)
can be proved exactly as in the case of Theorem \ref{ungm}. Instead of
Mahler's theorem \cite{mahler}, we now need to use that
$l_E(c) \not\in \bQ$ for any $0 \neq c \in pZ_{(p)}$, and the latter follows
by Bertrand's elliptic $p$-adic analogue of the Hermite-Lindemann
theorem \cite{ber}.
\qed

In the next result, we compute the group of solutions of the
$\delta_{\cP}$-character $\psi^{2e}_E$ given by (\ref{thechare}).

\begin{theorem}
Let $A$ be the $\delta_{\cP}$-ring $\bZ_{(\cP)}[\zeta_m]$ in
Example \ref{basicex}, and let $\psi^{2e}_{E,A}:E(A_{\cP}) \ra A_{\cP}$
be the homomorphism (\ref{indho}) induced by $\psi^{2e}_E$. Then
$$
{\rm Ker}\, \psi^{2e}_{E,A}=E(A_{\cP})_{tors}\, .
$$
\end{theorem}

This can be viewed as a version for several primes of the arithmetic analogue
\cite{frob} of Manin's Theorem of the Kernel \cite{manin}.

{\it Proof}. We use the same argument as in the proof of Theorem \ref{oastept}.
We use here the fact that the polynomials
$\phi_{p_l}^2-a_{p_l}\phi_{p_l}+p_l,\phi_{p_l}^M-1\in \bQ[\phi_{p_l}]$
are coprime. \qed

\section{Final Remarks and Questions}
\setcounter{theorem}{0}

We see an emerging pattern from the construction of $\delta_{\cP}$-characters
in the previous section. It can be best explained via the following:

\begin{definition}\label{equalizer}
Let $X$ be a quasi-projective smooth scheme over $\bZ_{(\cP)}$
with irreducible geometric fibers, and let $\cP=\{p_1,\ldots ,p_d\}$ be a
finite set of primes not dividing $N$.
Let $P$ be a uniform point of $X$. A $\delta_{\cP}$-function
$$
f =(f_{k}) \in \cO_{P,\cP}^r(X)
$$
will be said to have a {\it Dirac decomposition}
if there exists a family
$$
(\eta_{p_k}) \in \prod_{k=1}^d \cO(\cJ^s_{\{p_k\}}(X)\whk)\, ,
$$
$s \in \bZ_+$, and a family
$$
(\bar{\eta}_{\bar{p}_k}) \in \prod_{k=1}^d \bZ_{p_k}[\phi_l: \, l\in I_k]\, ,
$$
$I_k:=\{1,\ldots ,d\}\backslash \{p_k\}$, $\bar{p}_k=\cP\backslash p_k$,
such that $f_{k}=\bar{\eta}_{\bar{p}_k} \eta_{p_k}$ for all $k$.
\end{definition}

Then the following remarks are in order:

\begin{remark}
Any $\delta_{\cP}$-character of $\bG_a$, $\bG_m$, $E$ (where $E$ is an
elliptic curve over $\bZ_{(\cP)}$)
has a Dirac decomposition.
\qed
\end{remark}

\begin{remark}
\label{rdoi} Let $X$ be a projective curve over $\bZ_{(\cP)}$ of genus
$\geq 2$ with a uniform point $P$. Can we always find
$\delta_{\cP}$-functions on $X$ not belonging to $\bZ_{(\cP)}$? In
other words, do we have
$\cO_{P,\cP}^{\infty}(X)\neq \bZ_{(\cP)}$? (Note that by \cite{pjets} we
always have
$\cO(\cJ^1_{\{p_k\}}(X)\whk) \neq \bZ_{p_k}$ for each $k$.) Can we always
find a $\delta_{\cP}$-function in
$\cO_{P,\cP}^{\infty}(X)\backslash \bZ_{(\cP)}$ which has a Dirac
decomposition? By Theorem \ref{exE}, the answer to the last question
is positive if $X$ has a morphism to an elliptic curve over
$\bZ_{(\cP)}$; by Eichler-Shimura, this is always the case for
appropriate $\cP$s and  $X$ a modular curve of sufficiently high level.
\qed
\end{remark}

\begin{remark}
Let $X={\rm Spec}\, \bZ[1/6][a_4,a_6,(4a_4^3+27a_6^2)^{-1}]$. The theory of
differential modular forms
in the ODE setting \cite{difmod, book} is, essentially, the study of
some remarkable elements of the rings $\cO(\cJ^n_{\{p_k\}}(X)\whk)$
satisfying certain homogeneity properties. We may ask, in the arithmetic
PDE setting, for the existence and structure of {\it $\delta_{\cP}$-modular
forms} defined  as $\delta_{\cP}$-functions on $X$ with appropriate
homogeneity properties. Some such forms can be constructed via
Dirac decompositions, using Remark \ref{rdoi} above.
\qed
\end{remark}

\begin{remark}
 The construction  of non-trivial $\delta_{\cP}$-characters on
one dimensional groups via Dirac decompositions is a consequence
of the fact that the  ``ring of symbols,'' $\bZ_{(\cP)}[\phi_{\cP}]$, is
a commutative integral domain. Indeed, the key property used to
construct Dirac decompositions  is that

\begin{center}
\parbox{.1in}{(*)} \hspace{.2in}
\parbox{4in}{For any collection of non-zero elements
$\Theta_1,\ldots,\Theta_d$ in $\bZ_{(\cP)}[\phi_{\cP}]$, there exist
elements $\Lambda_1,\ldots ,\Lambda_d$ in $\bZ_{(\cP)}[\phi_{\cP}]$ such that
$\Lambda_1 \Theta_1=\cdots =\Lambda_d \Theta_d\neq 0$.}
\end{center}

This property is trivially true in any commutative integral domain.

Now, we can ask for a generalization of the theory of the present paper
in two possible directions:
\begin{enumerate}
\item[(a)] arithmetic PDEs in $0+d$ independent variables ($d \geq 2$),
 and $d_3\geq 2$ dependent variables (corresponding to commutative
algebraic groups of dimension $\geq 2$);
\item[(b)] arithmetic PDEs in $1+d$ independent variables
($d\geq 2$) and $d_3=1$ dependent variables (corresponding
to $\bG_a$, $\bG_m$, $E$).
\end{enumerate}
\medskip

However, in both cases, (a) and (b), the ``rings of symbols'' are
non-commutative and, in case (a), the corresponding ring has zero divisors.
Indeed, in case (a) the ``ring of symbols'' turns out to be  the ring of
$d_3\times d_3$ matrices with coefficients in $\bZ_{(\cP)}[\phi_{\cP}]$.
In case (b), the ``ring of symbols'' turns out to be the ring
$\bZ_{(\cP)}[\delta_q,\phi_{\cP}]$ generated by variables
$\delta_q,\phi_{p_1},...,\phi_{p_d}$ with the $\phi_{p_k}$s commuting
with each other, and with
$\delta_q \phi_{p_k}=p_k \cdot \phi_{p_k}\delta_q$ for all $k$;
cf. \cite{pde} for the case $d=1$.

In both cases, (a) and (b), property (*) fails for the corresponding rings
of symbols (and for the elements in those rings that are relevant to
our problem) so  the construction of characters with Dirac decompositions
fails as it is. An interesting problem would then be to determine all
the characters in cases (a) and (b)
respectively.
\qed
\end{remark}

\bibliographystyle{amsalpha}

\end{document}